\documentclass[12pt]{article}
\usepackage{amsmath, latexsym, amsfonts, amssymb, amsthm, amscd}
\usepackage[left=2.4cm, right=2cm, top=2.5cm, bottom=2.5cm]{geometry}
\usepackage{upquote}
\usepackage{color}
\usepackage{setspace}
\newtheorem{corollary}{Corollary}
\newtheorem{proposition}{Proposition}
\newtheorem{lemma}{Lemma}

\newtheorem{theorem}{Theorem}

\begin{document}

\title{Yule processes with rare mutation and their applications to percolation on $b$-ary trees.}
\author{G. Berzunza\footnote{ {\sc Institut f\"ur Mathematik, Universit\"at Z\"urich, Winterthurerstrasse 190, CH-8057 Z\"urich, Switzerland;} e-mail: gabriel.berzunza@math.uzh.ch}}
\maketitle

\vspace{0.2in}

\begin{abstract} 
\noindent We consider supercritical Bernoulli bond percolation on a large $b$-ary tree, 
in the sense that with high probability, there exists a giant cluster. We show that the size of the giant cluster
has non-gaussian fluctuations, which extends a result due to Schweinsberg \cite{Sch} in the case of random recursive trees. 
Using ideas in the recent work of Bertoin and Uribe Bravo \cite{Be1}, the approach developed in this work
relies on the analysis of the sub-population with ancestral type in a system of branching processes with
rare mutations, which may be of independent interest. This also allows us to establish the analogous result for
scale-free trees.  

\bigskip 


\end{abstract}

\section{Introduction and main result}

Consider a tree of large but finite size $n$ and perform Bernoulli bond percolation with parameter 
$p_{n} \in (0,1)$ that depends on the size of that tree. So each edge is removed with probability
$1-p_{n}$ and independently of the other edges, inducing a partition of the set of vertices into 
connected clusters. We are interested in the supercritical regime, in the sense that with high probability, 
there exists a giant cluster of size comparable to $n$, and its complement has also a size of order $n$.
In fact, it has been shown recently in \cite{Be2} that for fairly general families of trees,
the supercritical regime corresponds to parameters of the form $p_{n} = 1 - c/  \ell(n)$,
where $\ell(n)$ is an estimate of the height of a typical vertex in the structure. \\

In the case of the uniform random recursive trees (i.e. trees on an ordered set of vertices where the
smallest vertex serves as the root, and the sequence of vertices along any branch from
the root to a leaf is increasing) it easily seen that $\ell(n) = \ln n$, so choosing the 
percolation parameter so that
\begin{eqnarray} \label{ec1}
p_{n} =  1 - \frac{c}{\ln n},
\end{eqnarray}

\noindent where $c > 0$ is fixed, corresponds to the supercritical regime. More precisely, if $\Gamma_{n}$ denotes the size of the largest cluster, then
$\lim_{n \rightarrow \infty} n^{-1} \Gamma_{n} = e^{-c}$ in probability. This can be viewed as the 
law of large numbers for the giant cluster, and it is then natural to investigate its fluctuations.
Schweinsberg \cite{Sch} (see also Bertoin \cite{Be4} for an alternative approach) has shown that in 
this particular case, the fluctuations are non-Gaussian. Specifically
\begin{eqnarray} \label{ec19}
(n^{-1} \Gamma_{n} - e^{-c}) \ln n - ce^{-c} \ln \ln n \Rightarrow -ce^{-c} (\mathcal{Z}+ \ln c),
\end{eqnarray}

\noindent where $\Rightarrow$ means convergence in law  as $n \rightarrow \infty$ and the variable $\mathcal{Z}$ has the continuous Luria-Delbr\"uck distribution, i.e. its characteristic function is given by
\begin{eqnarray*}
\mathbb{E}(e^{i \theta \mathcal{Z}}) = \exp \left(- \frac{\pi}{2} |\theta| - i \theta \ln |\theta| \right), \hspace*{5mm} \theta \in \mathbb{R}.
\end{eqnarray*}

The main purpose of this work is to investigate the case of large random $b$-ary recursive trees ($b \geq 2$). 
The process to build a $b$-ary recursive tree starts at $n=1$ from the tree $T_{1}$ with one internal vertex (which corresponds to the 
root) and $b$ external vertices. Then, we suppose that $T_{n}$ has been constructed for some $n \geq 1$ that is a tree with $n$ internal vertices
and $(b-1)n+1$ external ones (also called leaves). Then choose an external vertex uniformly at random and replace it by an internal vertex to which $b$ new leaves are attached. 
In the case $b=2$, the algorithm yields a so-called binary search tree (see for instance Mahmoud \cite{Ma}, Drmota \cite{Drmota}). We consider that the size of the tree is the number of internal vertices. \\

Then we perform Bernoulli bond percolation with parameter given by (\ref{ec1}) on a random $b$-ary recursive tree of size $n$, which corresponds precisely to the supercritical regime as the case of the random recursive trees. Roughly speaking, since the $b$-ary recursive trees have also logarithmic height, i.e. the height of typical vertex is approximately $\ell(n) = (b \ln n) / (b-1)$ (see Javanian and Vahidi-Asl \cite{java}), one can verify that percolation then produces a giant cluster  whose size $C_{0}^{(p)}$ (number of internal vertices) satisfies
\begin{eqnarray*}
 \lim_{n \rightarrow \infty} n^{-1}C_{0}^{(p)} = e^{-\frac{b}{b-1}c}\hspace*{6mm} \text{in probability}.
\end{eqnarray*}

We now state the central result of this work, which shows that the fluctuations of the giant cluster 
in the case of the $b$-ary recursive trees are also described by the continuous Luria-Delbr\"uck distribution.
We stress that this distribution was further observed in relation with a random algorithm for the isolation of the root,
in the context of uniform random recursive tree by Iksanov and M\"ohle \cite{Iks}, and 
for random binary search tree by Holmgren \cite{Hol}.

\begin{theorem} \label{teo1}
 Set $\beta = b/(b-1)$, and assume that the percolation parameter $p_{n}$ is given by (\ref{ec1}). Then as $n \rightarrow \infty$, there is the weak convergence
 \begin{eqnarray*}
  (n^{-1}C_{0}^{(p)}- e^{- \beta c}) \ln n - \beta c e^{-\beta c} \ln \ln n \Rightarrow  -  \beta c e^{-\beta c}  \mathcal{Z}_{c,\beta}
 \end{eqnarray*}
 \noindent where
\begin{eqnarray} \label{ec8}
\mathcal{Z}_{c,\beta} =   \mathcal{Z} - \kappa_{\beta} + \ln ( \beta c )  
\end{eqnarray}

\noindent with $\mathcal{Z}$ having the continuous Luria-Delbr\"uck distribution,
\begin{eqnarray} \label{ec9}
\kappa_{\beta} = 1-\frac{1}{\beta} + \frac{1}{\beta} \sum_{k=2}^{\infty} \frac{(\beta)_{k}}{k!} \frac{(-1)^{k}}{k-1},  
\end{eqnarray}

\noindent and $(x)_{k} = x(x-1) \cdots (x-k+1)$, for $k \in \mathbb{N}$ and $x \in \mathbb{R}$, is the Pochhammer function. In particular, for $b = 2$, i.e. for the binary search tree case, $\kappa_{2} = 1$.
\end{theorem}

It should be noted the close similarity with the result for uniform recursive 
trees. It is remarkable that the normalizing functions and the limit in 
Theorem \ref{teo1} only depend on the parameter $\beta = b/(b-1)$ through some constants. Observe that the left-hand side of (\ref{ec19}) is the same as in Theorem \ref{teo1} for $\beta=1$; however the expressions (\ref{ec8}) and (\ref{ec9}) are not defined for $\beta=1$ ! \\

The basic idea of Schweinsberg \cite{Sch}, for establishing the result (\ref{ec19}) for uniform recursive trees relies on the estimation of the rate of decrease of the number of blocks in the Bolthausen-Sznitman coalescent, 
using the construction due to Goldschmidt and Martin \cite{Go} of the latter in terms of uniform recursive trees. On the other hand, the alternative approach of Bertoin \cite{Be4} makes use on the
remarkable coupling due to Iksanov and M\"ohle \cite{Iks} connecting the Meir and Moon algorithm
for the isolation of the root, with a certain random walk in the domain of attraction of the 
completely asymmetric Cauchy process. These approaches depend crucially
on the {\em splitting property} (see Section 3.1 in Bertoin \cite{Be3}) which fails for
the $b$-ary recursive trees. We thus have to use a different argument, although some guiding lines are
similar to \cite{Be4}. \\

Essentially, we consider a continuous time version of the growth algorithm of the $b$-ary tree which
bears close relations to Yule processes. The connection between recursive trees and branching processes 
is well-known, we make reference to Chauvin, et. al. \cite{Chauvin} for the binary search trees and
Bertoin and Uribe Bravo \cite{Be1} for the case of scale-free trees. In this way, we adapt the recent strategy
of \cite{Be1}. Roughly speaking, incorporating percolation to the algorithm yields systems of
branching processes with mutations, where a mutation event corresponds to 
disconnecting a leaf from its parent, and simultaneously replacing it by an internal vertex to which
$b$ new leaves are attached. Each percolation cluster size can then be thought of as a 
sub-population with some given genetic type. Hence the problem is reduced to study the fluctuations of the size of the sub-population with 
the ancestral type, which corresponds to the number of internal vertices connected to the root cluster. \\

The work is organized as follows. In Section \ref{sec2}, we introduce the system of branching
processes with rare mutations. We investigate the fluctuations of the size of the 
sub-population with the ancestral type, when the total population goes to infinity and 
the mutation parameter $1-p_{n}$ satisfies (\ref{ec1}). Then in Section \ref{sec3}, we make the link
with percolation on $b$-ary recursive trees in order to prove Theorem \ref{teo1}. 
Finally, we briefly show in Section \ref{sec4} that the present approach also applies to study the fluctuations of the size of the giant cluster for percolation on scale-free trees.

\section{Yule process with rare mutations} \label{sec2}

The purpose of this section is to introduce a system of branching process with rare mutations, 
which is quite similar to the one considered in \cite{Be1}, although there are also some
key differences (in particular, death may occur causing the extinction of sub-population with
the ancestral type). Then we focus on estimating the size of the sub-population with the ancestral type, when the total population in the 
system grows and the mutation parameter depends of the size of the latter. \\

We consider a population in which each individual is either a clone (i.e. an individual
with the ancestral type) or a mutant with some genetic type. A clone individual lives 
for an exponential time of parameter $1$, and gives birth at its death to $b$ clones with 
probability $p \in (0,1)$, or $b$ mutants that share the same genetic type with 
probability $1-p$. A mutant individual lives for an exponential time of parameter $1$,
and gives birth at its death to $b$ children of the same genetic of its parent. More precisely, the 
evolution of the population system is described by the process $\mathbf{Z}^{(p)} = (\mathbf{Z}^{(p)}(t):t \geq 0)$, where 
\begin{eqnarray*}
 \mathbf{Z}^{(p)}(t) = (Z_{0}^{(p)}(t), Z_{1}(t), \dots), \hspace*{5mm} \text{for} \hspace*{3mm} t \geq 0,
\end{eqnarray*}

\noindent is a collection of nonnegative variables which represents the current size
of the sub-populations. At the initial time, the sub-populations $Z_{i}(0)$ 
of type $i \geq 1$ are zero, and $Z_{0}^{(p)}(0) = b$ which is the size of the ancestral 
population. Formally, we take $\mathbf{Z}^{(p)}$ to be a pure-jump Markov chain whose transitions are described as follows. When at state 
$\mathbf{z}=(z_{i}:i \geq 0)$, our process jumps to a state $\tilde{\mathbf{z}} = (\tilde{z}_{i}: i \geq 0)$ 
where $\tilde{z}_{j} = z_{j}$ for $j \neq k$ and $\tilde{z}_{k} = z_{k} + (b-1)$ at rate
\[
 \left\{ \begin{array}{lcl}
              pz_{0} & \mbox{  if } & k =0, \\
              z_{k}  & \mbox{  if } & k \neq 0. \\
              \end{array}
    \right.
\]

\noindent This corresponds to a reproduction event 
 in the sub-population with type $k$. Otherwise, the process jumps from $\mathbf{z}$ to 
$\hat{\mathbf{z}}= (\hat{z}_{i}: i \geq 0)$ at rate $(1-p)z_{0}$ where, if $k$ is the first index such that $z_{k} = 0$, then 
 $\hat{z}_{0} = z_{0} -1$, $\hat{z}_{k}=b$, and $\hat{z}_{j} = z_{j}$ for $j \neq 0, k$. This corresponds to a mutation event of the sub-population with
 the ancestral type.  \\

The process of the total size of the population in the system
\begin{eqnarray*}
 Z(t) = Z_{0}^{(p)}(t) + \sum_{i \geq 1} Z_{i}(t), \hspace*{5mm} t \geq 0,
\end{eqnarray*}

\noindent is distributed as a Yule process, where each individual lives for an exponential time of parameter 1 and gives birth 
at its death to $b$ children, which then evolve independently of one another according to 
the same dynamics as their parent, no matter the choice of $p$.
Clearly, the process of the size of the sub-population with the ancestral 
type $Z_{0}^{(p)}$ is a continuous time branching
process, with reproduction law given by the distribution of $b \epsilon_{p}$, where $\epsilon_{p}$ stands
for a Bernoulli random variable with parameter $p$. Moreover, if for $i \geq 1$, we write
\begin{eqnarray*} 
 a_{i}^{(p)} = \inf \{ t \geq 0: Z_{i}(t) > 0 \},
\end{eqnarray*}

\noindent for the birth time of the sub-population with type $i$, then each process
\begin{eqnarray*}
 (Z_{i}(t-a_{i}^{(p)}): t \geq a_{i}^{(p)})
\end{eqnarray*}

\noindent is a branching process with the same reproduction law as $Z$ starting from $b$. 
Indeed, the different populations present in the system (i.e., those with strictly positive sizes)
evolve independently of one another. The following statement is just a formal formulation of 
the previous observation which should be plain from the construction of $\mathbf{Z}^{(p)}$; it 
is essentially Lemma 1 in \cite{Be1}.

\begin{lemma} \label{lema6}
 The processes $(Z_{i}(t-a_{i}^{(p)}): t \geq a_{i}^{(p)})$ for $i \geq 1$ form a sequence of
 i.i.d. branching process with the same law as $Z$ and with starting value $b$. Further, this sequence
 is independent of that of the birth-times $(a_{i}^{(p)})_{i \geq 1}$ and the process $Z_{0}^{(p)}$ of
 the sub-population with ancestral type.
\end{lemma}

We are now ready to present the main result of this section.  
We henceforth assume that the parameter $p = p_{n}$ is given by (\ref{ec1}) and for
simplicity, we write $p$ rather than $p_{n}$, omitting the integer $n$ from the notation. We consider the time
\begin{eqnarray*}
 \tau(n) = \inf \{t \geq 0: Z(t) = (b-1) n + 1 \}, 
\end{eqnarray*}

\noindent when the total population has size $(b-1) n + 1$. The size of the sub-population with the ancestral type 
at this time is given by
\begin{eqnarray*}
G_{n} := Z_{0}^{(p)}(\tau(n)).
\end{eqnarray*}

\begin{theorem} \label{teo2}
 Set $\beta = b/(b-1)$. As $n \rightarrow \infty$, there is the weak convergence
 \begin{eqnarray*}
  \left(n^{-1}G_{n}-\frac{1}{\beta-1} e^{-\beta c} \right) \ln n - \frac{\beta}{\beta-1} c e^{-\beta c} \ln \ln n \Rightarrow  -  \frac{\beta}{\beta-1} c e^{-\beta c} \left( \mathcal{Z}_{c,\beta} +1- \frac{1}{\beta}  \right),
 \end{eqnarray*}
 
 \noindent where $\mathcal{Z}_{c,\beta}$ is the random variable defined in (\ref{ec8}).
\end{theorem}

We stress that this result also allows us to deduce the fluctuations of the number 
of mutants in the total population, since this quantity is given by $(b-1) n +1 - G_{n}$. \\

The rest of this section is devoted to the proof of Theorem \ref{teo2}. 
Our approach is similar to that in \cite{Be4}. Broadly speaking, 
we divide the study of the fluctuations in two well-defined phases. The crucial point is to obtain a 
precise estimate of the number $\Delta_{n}$ of mutants when the total population of the system 
attains the size $(b-1)\lfloor \ln^{4} n \rfloor+1$; this can be viewed as the germ of the fluctuations of $(b-1) n + 1-G_{n}$. Then, we resume the growth of the system from size 
$(b-1)\lfloor \ln^{4} n \rfloor+1$ to the size $(b-1) n + 1$ and observe that 
the sub-population with the ancestral type grows essentially regularly. We point out that
even though the study of these two phases plays a key role in \cite{Be4}, 
the tools developed here to deal with each phase are much different from those used there.

\subsection{The germ of fluctuations} \label{germ}

In this first phase, we observe the growth of the system of branching processes until the time 
\begin{eqnarray*}
 \tau(\ln^{4} n ) = \inf \{t \geq 0: Z(t) = (b-1)\lfloor \ln^{4} n \rfloor +1\}, 
\end{eqnarray*}

\noindent which is when the total size of the population reaches $(b-1)\lfloor \ln^{4} n \rfloor+1$, and our purpose in this section is to estimate precisely the number $\Delta_{n}$ of mutants in the total population at this time, that is
\begin{eqnarray*}
 \Delta_{n} = (b-1)\lfloor \ln^{4}n \rfloor+1 - Z_{0}^{(p)}(\tau( \ln^{4} n)).
\end{eqnarray*}

We stress that the threshold $(b-1)\ln^{4} n+1$ is somewhat arbitrary, and any power close to $4$
of $\ln n$ would work just as well. However, as is remarked by Bertoin in \cite{Be4}, it is crucial
to choose a threshold which is both sufficiently high so that fluctuations are already visible, 
and sufficiently low so that one can estimate the germ with the desired accuracy. \\

We start by setting down the key results that lead us to the main result of this section, in order to give an easier articulation of the argument. 
In this direction, it is convenient to introduce the number $\Delta_{0, n}$ of mutants at time
\begin{eqnarray*}
 \tau_{0}(\ln^{4} n ) = \inf \{t \geq 0: Z_{0}^{(p)}(t) = (b-1)\lfloor \ln^{4} n \rfloor+1 \}, 
\end{eqnarray*}

\noindent which is when the size of the sub-population with the ancestral 
type reaches $(b-1)\lfloor \ln^{4} n \rfloor+1$, i.e.
\begin{eqnarray*}
 \Delta_{0,n} = Z(\tau_{0}( \ln^{4} n))-(b-1)\lfloor \ln^{4}n \rfloor -1.
\end{eqnarray*}

This will be useful since the distribution of $\Delta_{0,n}$ is easier to estimate than
that of $\Delta_{n}$. Then, we establish the following limit theorem in law that 
relates the fluctuations of $\Delta_{0, n}$ with the continuous Luria-Delbr\"uck variable $\mathcal{Z}$.

\begin{proposition} \label{pro1}
As $n \rightarrow \infty$, there is the weak convergence
\begin{eqnarray*}
\frac{\Delta_{0,n}}{\ln^{3}n} -  3 \frac{\beta}{\beta -1} c\ln \ln n \Rightarrow \frac{\beta}{\beta -1} c \left( \mathcal{Z}_{c,\beta} + 1 - \frac{1}{\beta}  \right)
\end{eqnarray*}

\noindent where $\mathcal{Z}_{c,\beta}$ is the random variable defined in (\ref{ec8}).
\end{proposition}

As we are interested in estimate the number $\Delta_{n}$
of mutants in the total population at time $\tau(\ln^{4}n)$, and we know the behavior
of $\Delta_{0,n}$, the purpose of the next
lemma is to point out that  these two quantities are close enough when $n \rightarrow \infty$. 
We need to introduce the notation:
\begin{eqnarray*}
 A_{n} = B_{n} + o(f(n)) \hspace*{5mm} \text{in probability},
\end{eqnarray*}

\noindent where $A_{n}$ and $B_{n}$ are two sequences of random variables and $f: \mathbb{N} \rightarrow (0, \infty)$
is a function, to indicate that $|A_{n} - B_{n}|/f(n) \rightarrow 0$ in probability when $n \rightarrow \infty$.

\begin{lemma} \label{lema8}
 We have
 \begin{eqnarray*}
  \Delta_{n} = \Delta_{0,n} + o(\ln^{3}n) \hspace*{5mm} \text{in probability}.
 \end{eqnarray*}
\end{lemma}

It then follows from Proposition \ref{pro1} that $\Delta_{n}$ and $\Delta_{0,n}$ have the 
same asymptotic behavior. Specifically:

\begin{corollary} \label{cor2}
As $n \rightarrow \infty$, there is the weak convergence
\begin{eqnarray*}
\frac{\Delta_{n}}{\ln^{3}n} -  3 \frac{\beta}{\beta-1} c\ln \ln n \Rightarrow \frac{\beta}{\beta -1} c \left( \mathcal{Z}_{c,\beta} + 1 - \frac{1}{\beta}  \right)
\end{eqnarray*}

\noindent where $\mathcal{Z}_{c,\beta}$ is the random variable defined in (\ref{ec8}).
\end{corollary}

The above result will be sufficient for our purpose. We now prepare the ground for the proofs
of Proposition \ref{pro1} and Lemma \ref{lema8}. Recall that we wish to study the behavior of the number $\Delta_{0,n}$ of mutants at time
$\tau_{0}(\ln^{4} n )$, which is easier than that of $\Delta_{n}$, thanks to Lemma \ref{lema6}.
In words, at time $\tau_{0}(\ln^{4} n )$ there is an independence property between the mutant sub-populations, and 
the process that counts the number of mutation events, which allows us to express $\Delta_{0,n}$ as a random 
sum of independent Yule processes. Clearly, 
the above is not possible at time $\tau(\ln^{4} n)$ due to the lack of independence within the sub-populations.
Formally, we start by writing 
\begin{eqnarray*}
 M(t) = \max \{i \geq 1: Z_{i}(t) > 0 \} 
\end{eqnarray*}

\noindent for the number of mutations that have occurred before time $t \geq 0$. Lemma \ref{lema6} ensures that $M$ is independent of the processes 
$(Z_{i}(t-a_{i}^{(p)}): t \geq a_{i}^{(p)})$ for $i \geq 1$. In addition, we note that
the jump times of $M$ are in fact $a_{1}^{(p)} < a_{2}^{(p)} < \cdots$. This enables us to express
the total mutant population at time $t$ as,
\begin{eqnarray*}
 Z_{\text{m}}(t) = \sum_{i=1}^{M(t)} Z_{i}(t - a_{i}^{(p)}),
\end{eqnarray*}

\noindent and we are thus interested in
\begin{eqnarray} \label{ec17}
 \Delta_{0,n} = Z_{\text{m}}(\tau_{0}( \ln^{4} n)).
\end{eqnarray}

We now turn our attention to study the fluctuations of $\Delta_{0,n}$ through 
the analysis of its characteristic function. 
In this direction, we will be mainly interested in the following feature of $Z_{m}(t)$.

\begin{lemma} \label{lema7}
We have for $t \geq 0$ and $\theta \in \mathbb{R}$.
\begin{itemize} 
 \item[i)] The characteristic function of $Z(t)$ started from $Z(0) = b$,
\begin{eqnarray} \label{ec11}
 \varphi_{t}(\theta) = \mathbb{E}\left[ e^{i \theta Z(t)} \big| Z(0) = b \right] = \left( \frac{e^{i \theta(b-1)} e^{-(b-1)t}}{  1 - e^{i \theta(b-1)} + e^{i \theta(b-1)} e^{-(b-1)t}} \right)^{\frac{b}{b-1}}.
\end{eqnarray} 
 \item[ii)] We have
  \begin{eqnarray} \label{ec10}
  \mathbb{E}[e^{i \theta Z_{\text{m}}(t)}] = \mathbb{E} \left[ \exp \left( (1-p) \int_{0}^{t} Z_{0}^{(p)}(t-s) \left( \varphi_{s}(\theta) - 1 \right) {\rm d} s  \right) \right].
\end{eqnarray}
\end{itemize}
\end{lemma}

\noindent {\bf Proof:} Recall that the processes $(Z_{i}(t - a_{i}^{(p)}): t \geq a_{i}^{(p)})$ for $i \geq 1$ are 
i.i.d. branching process with the same law as $Z$ with starting value $b$. Then according to 
page 109 in Chapter III of Athreya and Ney \cite{Athe}, their characteristic function
is given by the expression (\ref{ec11}). We now observe from the dynamics of $\mathbf{Z}^{(p)}$ that the counting process
$M$ has jumps at rate $(1-p)Z_{0}^{(p)}$. Moreover, conditionally on $Z_{0}^{(p)}$,
the process $Z_{\text{m}}$ is a non homogeneous filtered Poisson process whose
characteristic function can be written in terms of the characteristic function of $Z_{i}$. By
extending equation (5.43) of Parzen \cite{parzen} slightly to allow the underlying Poisson process to be non homogeneous, we obtain 
 \begin{eqnarray} \label{e1}
  \mathbb{E} \left[e^{i \theta Z_{\text{m}}(t)} \big| \left(Z_{0}^{(p)}(s):0 \leq s \leq t \right) \right] =  \exp \left( (1-p) \int_{0}^{t} Z_{0}^{(p)}(s) \left( \varphi_{t-s}(\theta) - 1 \right) {\rm d} s  \right), 
\end{eqnarray}

\noindent for $t\geq 0$ and $\theta \in \mathbb{R}$. Hence our claim follows after taking expectations on both sides of the equation and make a simple change of variables.  \hfill $\Box$ \\

We recall some important properties of the branching processes $Z$ and $Z_{0}^{(p)}$, 
which will be useful later on. The process
\begin{eqnarray*}
 W(t) := e^{-(b-1)t}Z(t), \hspace*{6mm} t \geq 0
\end{eqnarray*}

\noindent is a nonnegative square-integrable martingale which converges a.s. and in
$L^{2}(\mathbb{P})$, and we write $W(\infty)$ for its terminal value.
Furthermore $W(\infty) > 0$ a.s. since $Z$ can not become extinct (we also pointed out that $Z$ never explodes a.s.). 
Similarly, the process
\begin{eqnarray*}
 W_{0}^{(p)}(t) = e^{-(bp-1)t}Z_{0}^{(p)}(t), \hspace*{6mm} t \geq 0
\end{eqnarray*}

\noindent is a martingale which terminal value is denoted by $W_{0}^{(p)}(\infty)$. In addition, following the proof of Lemma 3 in \cite{Be1} we have  
\begin{lemma} \label{lema4}
 It holds that 
 \begin{eqnarray*}
  \lim_{p \rightarrow 1, t \rightarrow \infty} \mathbb{E}_{z} \left[ \sup_{s \geq t} \left| W_{0}^{(p)}(s) - W(\infty) \right|^{2} \right] = 0.
 \end{eqnarray*}
 \noindent In particular, $W_{0}^{(p)}(\infty)$ converges to $W(\infty)$ in $L^{2}(\mathbb{P})$ as $p \rightarrow 1$.
\end{lemma}

We next estimate the characteristic function of $Z_{\text{m}}(t)$ given in 
(\ref{ec10}), but we still need some additional notation. For $t \geq 0$,
\begin{eqnarray*}
 I^{(p)}(t) = (1-p)\int^{t}_{0} Z_{0}^{(p)}(t-s)(\varphi_{s}(u)-1 ) {\rm d} s
\end{eqnarray*}

\noindent and
\begin{eqnarray*}
 I^{(p)}_{\text{m}}(t) = (1-p) W_{0}^{(p)}(\infty) e^{(b-1)t} \int^{t}_{0} e^{-(b-1)s} (\varphi_{s}(u)-1 ) {\rm d} s,
\end{eqnarray*} 

\noindent where $u = \theta / \left(\beta c \ln^{3} n \right) $ for $\theta \in \mathbb{R}$ and $\beta = b /(b-1)$.

\begin{lemma} \label{lema1}
We have
\begin{eqnarray*}
 \lim_{n \rightarrow \infty} \left(  I^{(p)}(\tau_{0}(\ln^{4} n )) -  I^{(p)}_{{\rm m}}(\tau_{0}(\ln^{4} n )) \right) = 0  \hspace*{6mm} \text{in probability}.
 \end{eqnarray*}
\end{lemma}

\noindent {\bf Proof:} Define the function 
\begin{eqnarray*}
 I^{(p)}_{\text{a}}(t) = (1-p) W_{0}^{(p)}(\infty)e^{(bp-1)t}  \int_{0}^{t} e^{-(bp-1)s} (\varphi_{s}(u)-1 ) {\rm d} s, \hspace*{6mm} t \geq 0,
\end{eqnarray*} 

\noindent which is simply obtained by replacing $b$ by $bp$ in the exponential terms of $I^{(p)}_{\text{m}}(t)$. 
We first prove that 
\begin{eqnarray} \label{ec16}
 \lim_{n \rightarrow \infty}  \left( I^{(p)}(\tau_{0}(\ln^{4} n)) -  I^{(p)}_{\text{a}}(\tau_{0}(\ln^{4} n)) \right)  = 0 \hspace*{5mm} \text{in probability}.
\end{eqnarray}

In this direction, we observe from the triangle inequality that
\begin{align} \label{ec12}
 & \left| I^{(p)}( \tau_{0}(\ln^{4} n ) ) -  I^{(p)}_{\text{a}}( \tau_{0}(\ln^{4} n ) )  \right| ~~~~~~ \nonumber \\
  &~~~ \leq   (1-p) e^{(bp-1) \tau_{0}(\ln^{4} n )} \int_{0}^{ \tau_{0}(\ln^{4} n ) } | W_{0}^{(p)}(\tau_{0}(\ln^{4} n ) -s)-W_{0}^{(p)}(\infty)||\varphi_{s}(u)-1| e^{-(bp-1)s} {\rm d}s.  
\end{align} 

\noindent We define
\begin{eqnarray*}
 A^{(p)} :=  \frac{3}{2(bp-1)}\sup_{s \geq 0} e^{(bp-1)s/3} |W_{0}^{(p)}(s)-W_{0}^{(p)}(\infty)|,
\end{eqnarray*}

\noindent and since Lemma 2 in \cite{Be1} shows that $A^{(p)}$
is bounded in $L^{2}(\mathbb{P})$, we have by the Markov inequality that
\begin{eqnarray}\label{ec15}
 \lim_{n \rightarrow \infty} \left( \ln^{-\frac{1}{3}}n \right)  A^{(p)} = 0 \hspace*{6mm} \text{in probability}.
\end{eqnarray}

We set $t_{n} =  (b-1)^{-1} \ln \ln n$ and recall that $\varphi_{t}(u)$ fulfills (\ref{ec11}). Hence from the inequality $|e^{ix}-1| \leq 2$ for $x \in \mathbb{R}$, we have that 
\begin{eqnarray} \label{rec1}
|\varphi_{t}(u)-1| \leq 2.
\end{eqnarray}
 Then, 
\begin{align} \label{ec13}
  & (1-p) e^{(bp-1)\tau_{0}(\ln^{4} n )} \int_{\tau_{0}(\ln^{4} n ) - t_{n}}^{\tau_{0}(\ln^{4} n )} | W_{0}^{(p)}(\tau_{0}(\ln^{4} n ) - s)-W_{0}^{(p)}(\infty)||\varphi_{s}(u)-1| e^{-(bp-1)s} {\rm d}s  \nonumber \\
  &~~~~~~~~~~~~~~~~~~~~~~~~~~~~~~~~~~ \leq 2(1-p)  \int_{0}^{t_{n}} | W_{0}^{(p)}(s)-W_{0}^{(p)}(\infty)| e^{(bp-1)s} {\rm d}s \nonumber \\
  &~~~~~~~~~~~~~~~~~~~~~~~~~~~~~~~~~~ \leq 2(1-p) \left(\ln^{\frac{2}{3}} n \right) A^{(p)}. 
\end{align}

\noindent On the other hand, from (\ref{ec11}) and the inequality $|e^{ix}-1| \leq |x|$ for $x \in \mathbb{R}$, we get that 
\begin{eqnarray} \label{rec2}
 |\varphi_{t}(u)-1| \leq b|u|e^{(b-1)t},
\end{eqnarray}

\noindent which implies that
\begin{align} \label{ec14}
  & (1-p) e^{(bp-1)\tau_{0}(\ln^{4} n )} \int_{0}^{\tau_{0}(\ln^{4} n ) - t_{n}} | W_{0}^{(p)}(\tau_{0}(\ln^{4} n ) - s)-W_{0}^{(p)}(\infty)||\varphi_{s}(u)-1| e^{-(bp-1)s} {\rm d}s \nonumber \\
  &~~~~~~~~~~~~~~~~ \leq  (1-p)b|u| e^{(b-1)\tau_{0}(\ln^{4} n )} \int_{t_{n}}^{\tau_{0}(\ln^{4} n)} | W_{0}^{(p)}(s)-W_{0}^{(p)}(\infty)| e^{-(b-1)s} e^{(bp-1)s} {\rm d}s \nonumber \\
  &~~~~~~~~~~~~~~~~ \leq   2b|u|(1-p)\left(\ln^{-\frac{1}{3}} n \right)  A^{(p)} e^{(b-1)\tau_{0}(\ln^{4} n )}. 
\end{align} 

We recall that $u = \theta / \left( \beta c \ln^{3} n \right)$, and $p=p_{n}$ is given by (\ref{ec1}). Then from (\ref{ec12}), (\ref{ec13}), and (\ref{ec14})  follow that
\begin{eqnarray*} 
    \left| I^{(p)}( \tau_{0}(\ln^{4} n ) ) -  I^{(p)}_{\text{a}}( \tau_{0}(\ln^{4} n ) )  \right|  & \leq &  2\left(  c   + (b-1)|\theta| \left( \ln^{-4}n \right)  e^{(b-1)\tau_{0}(\ln^{4} n )} \right) \left(\ln^{-\frac{1}{3}}n \right)A^{(p)}.
\end{eqnarray*} 

\noindent We observe that since $Z_{0}^{(p)}(\tau_{0}(\ln^{4} n )) = (b-1)\lfloor \ln^{4}n \rfloor+1$, Lemma  \ref{lema4} ensures that
\begin{eqnarray} \label{rec4}
 \lim_{n \rightarrow \infty} \left( \ln^{-4}n \right)  e^{(b-1)\tau_{0}(\ln^{4} n )} = \frac{b-1}{W(\infty)} \hspace*{6mm} \text{in probability},
\end{eqnarray}

\noindent where $W(\infty)$ is strictly positive almost surely.  Thus, we deduce (\ref{ec16}) from (\ref{ec15}) by letting $n \rightarrow \infty$ in the last inequality. \\

 Next, we prove that
\begin{eqnarray} \label{rec5}
 \lim_{n \rightarrow \infty}   \left( I^{(p)}_{\text{m}}(\tau_{0}(\ln^{4} n )) - I^{(p)}_{\text{a}}(\tau_{0}(\ln^{4} n ))\right)  = 0 \hspace*{6mm} \text{in probability},
\end{eqnarray}

\noindent by proceeding similarly as the proof of (\ref{ec16}). We observe for the triangle inequality that
\begin{align} \label{rrec1}
 & \left| I^{(p)}_{\text{m}}( \tau_{0}(\ln^{4} n ) ) -  I^{(p)}_{\text{a}}( \tau_{0}(\ln^{4} n ) )  \right| ~~~~~~ \nonumber \\
  &~~~~~~~~ \leq   (1-p) W_{0}^{(p)}(\infty) \int_{0}^{ \tau_{0}(\ln^{4} n ) } 
   \left | 1- e^{b(1-p)(s-\tau_{0}(\ln^{4} n ))} \right| |\varphi_{s}(u)-1| e^{-(b-1)(s- \tau_{0}(\ln^{4} n ))} {\rm d}s.  
\end{align} 

\noindent We set $t_{n} =  (b-1)^{-1} \ln \ln n$ again. Hence from the inequality (\ref{rec1}) we have that
\begin{align*} 
  & (1-p)W_{0}^{(p)}(\infty)  \int_{\tau_{0}(\ln^{4} n ) - t_{n}}^{\tau_{0}(\ln^{4} n )} \left | 1- e^{b(1-p)(s-\tau_{0}(\ln^{4} n ))} \right| |\varphi_{s}(u)-1| e^{-(b-1)(s- \tau_{0}(\ln^{4} n ))}  {\rm d}s  \nonumber \\
  &~~~~~~~~~~~ \leq 2(1-p)W_{0}^{(p)}(\infty)  \int_{\tau_{0}(\ln^{4} n ) - t_{n}}^{\tau_{0}(\ln^{4} n )} \left | 1- e^{b(1-p)(s-\tau_{0}(\ln^{4} n ))} \right|  e^{-(b-1)(s- \tau_{0}(\ln^{4} n ))}  {\rm d}s  \nonumber  \\
  &~~~~~~~~~~~ = 2(1-p)W_{0}^{(p)}(\infty)  \int_{0}^{t_{n}} \left | 1- e^{-b(1-p)s} \right|  e^{(b-1)s}  {\rm d}s.
\end{align*}

\noindent Then by making the change of variables $x = e^{(b-1)s}$, we get that
\begin{align} \label{rec3}
  & (1-p)W_{0}^{(p)}(\infty)  \int_{\tau_{0}(\ln^{4} n ) - t_{n}}^{\tau_{0}(\ln^{4} n )} \left | 1- e^{b(1-p)(s-\tau_{0}(\ln^{4} n ))} \right| |\varphi_{s}(u)-1| e^{-(b-1)(s- \tau_{0}(\ln^{4} n ))}   {\rm d}s  \nonumber \\
  &~~~~~~~~~~~ \leq  \frac{2}{b-1}(1-p)W_{0}^{(p)}(\infty)  \int_{1}^{\ln n} \left (1- x^{-\beta(1-p)} \right)   {\rm d}x \nonumber \\
  &~~~~~~~~~~~ \leq  \frac{2}{b-1}(1-p)W_{0}^{(p)}(\infty) \left( 1- (\ln n)^{-\beta(1-p)}\right)\ln n. 
\end{align}

\noindent On the other hand, from the inequality (\ref{rec2}) we have that 
\begin{align} \label{rrec14}
  & (1-p) W_{0}^{(p)}(\infty) \int_{0}^{\tau_{0}(\ln^{4} n ) - t_{n}}  \left | 1- e^{b(1-p)(s-\tau_{0}(\ln^{4} n ))} \right| |\varphi_{s}(u)-1| e^{-(b-1)(s- \tau_{0}(\ln^{4} n ))}  {\rm d}s \nonumber \\
  &~~~~~~~~~~~~~~~~ \leq  (1-p)b|u| e^{(b-1)\tau_{0}(\ln^{4} n )} W_{0}^{(p)}(\infty) \int_{t_{n}}^{\tau_{0}(\ln^{4} n)}  \left ( 1- e^{-b(1-p)s} \right) {\rm d} s \nonumber \\ 
  &~~~~~~~~~~~~~~~~ \leq \frac{b^{2}}{2} (1-p)^{2} |u| e^{(b-1)\tau_{0}(\ln^{4} n )} W_{0}^{(p)}(\infty) (\tau_{0}(\ln^{4} n ))^{2}.
  \end{align} 
  
\noindent Recall that $u = \theta / \left( \beta c \ln^{3} n \right)$, and $p=p_{n}$ is given by (\ref{ec1}). Then from (\ref{rrec1}), (\ref{rec3}), and (\ref{rrec14})  follow that
\begin{align*} 
  &  \left| I^{(p)}_{\text{m}}( \tau_{0}(\ln^{4} n ) ) -  I^{(p)}_{\text{a}}( \tau_{0}(\ln^{4} n ) )  \right|  \nonumber \\
&~~~~~~~~~~~ \leq 2b(b-1)c W_{0}^{(p)}(\infty) \left(  1- e^{- \beta c \frac{\ln \ln n}{\ln n} } +  |\theta| (\tau_{0}(\ln^{4} n ))^{2} e^{(b-1)\tau_{0}(\ln^{4} n )} \ln^{-5}n \right). 
  \end{align*} 

We deduce from (\ref{rec4}) that
\begin{eqnarray*} 
\lim_{n \rightarrow \infty} \frac{\tau_{0}(\ln^{4} n )}{4(b-1)^{-1}\ln \ln n} =1 \hspace*{6mm} \text{in probability},
\end{eqnarray*}

\noindent and since $\lim_{n \rightarrow \infty} W^{(p)}_{0}(\infty) = W(\infty)$ in $L^{2}(\mathbb{P})$ from Lemma \ref{lema4}, we get (\ref{rec5}) from (\ref{rec4}) by letting $n \rightarrow \infty$ in the last inequality. Finally, our claim follows by combining (\ref{ec16}) and (\ref{rec5}). \hfill $\Box$ \\

We observe that thanks to (\ref{ec11}), the integral $I_{\text{m}}^{(p)}$ can be computed explicitly.

\begin{lemma} \label{lema5}
 We have for $t \geq 0$,
 \begin{eqnarray*}
 \int^{t}_{0} e^{-(b-1)s} (\varphi_{s}(u)-1 ) {\rm d} s = \frac{1-e^{iu(b-1)}}{(b-1)e^{iu(b-1)}} \left( \beta \ln(1 - e^{iu(b-1)} + e^{iu(b-1)} e^{-(b-1)t}) + \kappa_{b,u}(t)\right),
\end{eqnarray*} 

\noindent where 
\begin{eqnarray} \label{ec18}
 \kappa_{b,u}(t) = \sum_{k=2}^{\infty} \frac{\left(\beta\right)_{k}}{k!}   \frac{(e^{iu(b-1)}-1)^{k-1}}{k-1} \left( 1 - \frac{1}{(1-e^{iu(b-1)}+e^{iu(b-1)}e^{-(b-1)t})^{k-1}} \right),
\end{eqnarray}

\noindent and $(\cdot)_{k}$ is the Pochhammer function.
\end{lemma}

\noindent {\bf Proof:} Define the function
\begin{eqnarray*}
 f(\lambda) = \int^{t}_{0} e^{-(b-1)r} \left(\left( \frac{e^{-\lambda (b-1)} e^{-(b-1)r}}{  1 - e^{-\lambda (b-1)} + e^{-\lambda (b-1)} e^{-(b-1)r}} \right)^{\beta}-1 \right) {\rm d} r, \hspace*{4mm} \lambda \geq 0.
\end{eqnarray*}

Hence by setting $x =  1 - e^{-\lambda (b-1)} + e^{-\lambda (b-1)} e^{-(b-1)r}$ and $y_{\lambda} = e^{-\lambda (b-1)}$ for convenience, we have that
\begin{eqnarray*}
 f(\lambda) = \frac{1}{(b-1)y_{\lambda}}\int^{1}_{1 - y_{\lambda} + y_{\lambda} e^{-(b-1)t}}  \left(\left( \frac{x+y_{\lambda}-1}{  x} \right)^{\beta}-1 \right) {\rm d}  x.
\end{eqnarray*}

\noindent Moreover, using a well-known extension of Newton's binomial formula, we get
\begin{eqnarray*}
 f(\lambda) = \frac{1}{(b-1)y_{\lambda}} \sum_{k=1}^{\infty} \frac{(\beta)_{k}}{k!}(y_{\lambda}-1)^{k} \int^{1}_{1 - y_{\lambda} + y_{\lambda} e^{-(b-1)t}} x^{-k} {\rm d} x,
\end{eqnarray*}

\noindent where $(\cdot)_{k}$ is the Pochhammer function. Note that the series converges
absolutely since $\beta > 0$ and $|y_{\lambda} - 1|/x \leq 1$,
 for $1 - y_{\lambda} + y_{\lambda} e^{-(b-1)t} \leq x \leq 1$. 
 Then straightforward calculations yield
\begin{eqnarray*}
 f(\lambda) = \frac{1-y_{\lambda}}{(b-1)y_{\lambda}} \left( \beta \ln(1 - y_{\lambda} + y_{\lambda} e^{-(b-1)t}) + \kappa_{b,\lambda}^{\prime}(t)\right),
\end{eqnarray*}

\noindent where 
\begin{eqnarray*}
 \kappa_{b,\lambda}^{\prime}(t) = \sum_{k=2}^{\infty} \frac{\left(\beta\right)_{k}}{k!}   \frac{(y_{\lambda}-1)^{k-1}}{k-1} \left( 1 - \frac{1}{(1-y_{\lambda}+y_{\lambda}e^{-(b-1)t})^{k-1}} \right).   
\end{eqnarray*}

We note that the function $f$ allows an analytic extension to $\{\lambda \in \mathbb{C}: \text{Re} \hspace*{1mm} \lambda \geq 0 \}$. 
Then by taking into account the principal branch of the complex logarithm, we conclude that 
\begin{eqnarray*}
 f(-i u) = \int^{t}_{0} e^{-(b-1)r} \left(\left( \frac{e^{i u(b-1)} e^{-(b-1)r}}{  1 - e^{i u(b-1)} + e^{i u(b-1)} e^{-(b-1)r}} \right)^{\beta}-1 \right) {\rm d} r,
\end{eqnarray*}

\noindent and our assertion follows by observing that $\kappa_{b,\lambda}^{\prime}(t) = \kappa_{b,u}(t)$ when $\lambda = -iu$.  \hfill $\Box$ \\

We are now able to establish Proposition \ref{pro1}. \\

\noindent {\bf Proof of  Proposition \ref{pro1}:} Fix $\theta \in \mathbb{R}$ and define $m_{n} = \beta c \ln^{3}n$. 
From the identity (\ref{ec17}) and Lemma \ref{lema7}, the characteristic function of 
$m_{n}^{-1} \Delta_{0,n} -(\beta-1)^{-1} \ln m_{n}$ is given by
\begin{eqnarray*}
\mathbb{E}\left[e^{i \theta (m_{n}^{-1} \Delta_{0,n} - (\beta-1)^{-1}\ln m_{n})}\right] & = & \mathbb{E}\left[  \exp\left( I^{(p)}(\tau_{0}(\ln^{4} n )) -i\theta(\beta-1)^{-1} \ln m_{n}  \right) \right].
\end{eqnarray*}

\noindent Recall that by Lemma \ref{lema1} we have
\begin{eqnarray*}
 \lim_{n \rightarrow \infty} \left( I^{(p)}(\tau_{0}(\ln^{4} n )) -  I^{(p)}_{\text{m}}(\tau_{0}(\ln^{4} n )) \right)= 0 \hspace*{5mm} \text{in probability.}
 \end{eqnarray*}

\noindent Then, we must verify that
\begin{align} \label{ec3}
  & \lim_{n \rightarrow \infty} \left(  I^{(p)}_{\text{m}} (\tau_{0}(\ln^{4} n )) - i\theta (\beta-1)^{-1} \ln m_{n} \right)  \nonumber \\
  &~~~~~~ = -  i\theta (\beta-1)^{-1} \left( \kappa_{\beta} - 1 + \frac{1}{\beta} \right)  - i \theta (\beta-1)^{-1} \ln |(\beta-1)^{-1}\theta| - \frac{1}{2} \pi |(\beta-1)^{-1}\theta| 
 \end{align} 

\noindent in probability. In this direction, we define $y_{u} = e^{iu(b-1)}$ for convenience and recall also from Lemma \ref{lema5} that
\begin{align*}
  & I^{(p)}_{\text{m}}(\tau_{0}(\ln^{4} n )) \nonumber \\
  &~~ = (1-p) W_{0}^{(p)}(\infty) e^{(b-1)\tau_{0}(\ln^{4} n )} \frac{1-y_{u}}{(b-1)y_{u}} \left( \beta \ln(1 - y_{u} + y_{u} e^{-(b-1)\tau_{0}(\ln^{4} n )}) + \kappa_{b,u}(\tau_{0}(\ln^{4} n ))\right),
 \end{align*} 

\noindent  where $\kappa_{b,u}(\tau_{0}(\ln^{4} n ))$ is defined in (\ref{ec18}) and $u = \theta / \left( \beta c \ln^{3} n \right)$. We know from Lemma \ref{lema4} that
\begin{eqnarray*}
 \lim_{n \rightarrow \infty} e^{-(b-1)\tau_{0}(\ln^{4} n )} Z(\tau_{0}(\ln^{4} n )) = \lim_{n \rightarrow \infty} e^{-(bp-1)\tau_{0}(\ln^{4} n )} Z_{0}^{(p)}(\tau_{0}(\ln^{4} n )) = W(\infty) 
\end{eqnarray*}

\noindent in probability, and $\lim_{n \rightarrow \infty} W^{(p)}_{0}(\infty) = W(\infty)$ in $L^{2}(\mathbb{P})$. Hence since $p=p_{n}$ fulfilled (\ref{ec1}), we deduce that 
\begin{eqnarray*}
\lim_{n \rightarrow \infty} \frac{(1-p) W_{0}^{(p)}(\infty) e^{(b-1)\tau_{0}(\ln^{4} n )}}{ (b-1)c\ln^{3}n} =1 \hspace*{6mm} \text{in probability}.
\end{eqnarray*}

On the other hand,
\begin{eqnarray*}
y_{u} = 1 + O \left( \frac{1}{m_{n}} \right) \hspace*{5mm} \text{and} \hspace*{5mm} m_{n}\left(1-y_{u}\right) = -i\theta(b-1) + O \left(\frac{1}{m_{n}} \right),
\end{eqnarray*}

\noindent and since $b-1 = (\beta -1)^{-1}$, we deduce that
\begin{align*}
  & \lim_{n \rightarrow \infty} \left( I^{(p)}_{m}(\tau_{0}(\ln^{4} n )) -i\theta (\beta -1)^{-1} \ln m_{n}  \right)  \nonumber \\
  &~~~~~~~~~~~~~ = - i \theta (\beta -1)^{-1} \ln(-i\theta(\beta -1)^{-1}) -i \theta (\beta -1)^{-1} \left( \kappa_{\beta} - 1 + \frac{1}{\beta} \right) 
 \end{align*} 

\noindent in probability, which implies (\ref{ec3}). Finally, we observe from 
(\ref{e1}) and the modulus inequality for conditional expectation that 
\begin{eqnarray*}
\left | \exp\left( I^{(p)}(\tau_{0}(\ln^{4} n )) -i\theta(\beta-1)^{-1} \ln m_{n} \right) \right | \leq 1.
\end{eqnarray*}

\noindent Therefore, by the dominated convergence theorem we conclude that the Fourier transform of $m_{n}^{-1} \Delta^{0}_{n} - (\beta -1)^{-1}\ln m_{n}$ converges pointwise  as $n$ tends to infinity to the continuous function 
\begin{eqnarray*}
 \theta \mapsto \exp \left( -  i\theta(\beta -1)^{-1}\left( \kappa_{\beta} - 1 + \frac{1}{\beta} \right)  - i \theta (\beta -1)^{-1}\ln |(\beta -1)^{-1}\theta| - \frac{1}{2} \pi |(\beta -1)^{-1}\theta|\right),
\end{eqnarray*}

\noindent and then our claim follows for the continuity theorem for Fourier transforms. \hfill $\Box$ \\

We now turn our attention to the proof of Lemma \ref{lema8}. \\

We imagine that we begin our observation of the system of branching
processes with rare mutations $\mathbf{Z}^{(p)}$ once it has reached the size $(b-1)\lfloor \ln^{4}n \rfloor+1$, that is,
from the time $\tau(\ln^{4} n)$. We thus write $Z^{\prime} = (Z^{\prime}(t): t \geq 0)$ for the process
of the total size of the population started from $Z^{\prime}(0) = (b-1)\lfloor \ln^{4}n \rfloor+1$, which has the same
law of the Yule process $Z$ described at the beginning of Section \ref{sec2}. We introduce the time
\begin{eqnarray*}
 \tau^{\prime}( \ln^{4}n ) = \inf \{t \geq 0: Z^{\prime}(t) = Z( \tau_{0}( \ln^{4}n  ) )\},
\end{eqnarray*}

\noindent at which it hits $Z( \tau_{0}( \ln^{4}n  ))$. Equivalently, $\tau^{\prime}( \ln^{4}n )$ is the time needed to have a population with the ancestral type of size $(b-1)\lfloor \ln^{4}n \rfloor+1$.
We shall first estimate this quantity. 

\begin{lemma} \label{lema9}
We have
 \begin{eqnarray*}
 \lim_{ n \rightarrow \infty} \tau^{\prime}( \ln^{4}n ) = 0 \hspace*{6mm} \text{in probability}.
 \end{eqnarray*}
\end{lemma}

\noindent {\bf Proof:} We know that
\begin{eqnarray*}
 \lim_{n \rightarrow \infty} e^{-(b-1)\tau(\ln^{4}n)}Z(\tau(\ln^{4}n)) = W(\infty) \hspace*{5mm} \text{in probability}.
\end{eqnarray*}

\noindent By definition of the time $\tau(\ln^{4}n)$, we have $Z(\tau(\ln^{4}n)) =  (b-1)\lfloor \ln^{4}n \rfloor+1$, hence
\begin{eqnarray*}
\lim_{n \rightarrow \infty} \frac{\tau(\ln^{4}n)}{4(b-1)^{-1}\ln \ln n}=1 \hspace*{5mm} \text{in probability}.
\end{eqnarray*}

On the other hand, from Lemma \ref{lema4} we have that
\begin{eqnarray*}
 \lim_{n \rightarrow \infty} e^{-(bp-1)\tau_{0}(\ln^{4}n)}Z_{0}^{(p)}(\tau_{0}(\ln^{4}n)) = W(\infty) \hspace*{5mm} \text{in probability}.
\end{eqnarray*}

\noindent Recall that $p=p_{n}$ is given by (\ref{ec1}), and observe that $Z_{0}^{(p)}(\tau_{0}(\ln^{4}n)) = (b-1)\lfloor \ln^{4}n \rfloor+1$. Hence
\begin{eqnarray*}
 \lim_{n \rightarrow \infty} e^{-(bp-1) (\tau(\ln^{4}n) -\tau_{0}(\ln^{4}n))} = 1 \hspace*{5mm} \text{in probability},
\end{eqnarray*}

\noindent and our claim follows from the identity $\tau^{\prime}( \ln^{4}n ) = \tau_{0}(\ln^{4}n) - \tau(\ln^{4}n)$. \hfill $\Box$ \\

We observe that the population at time $\tau(\ln^{4} n)$ when we start our observation consists of $\Delta_{n}$ mutants and $(b-1)\lfloor \ln^{4} n \rfloor+1 - \Delta_{n}$
individuals of the ancestral type. Then, we write $Z^{\prime (p)}_{0} = (Z^{\prime (p)}_{0}(t): t \geq 0)$ for the process
that counts the number of individuals with the ancestral type, which has the same law as $Z_{0}^{(p)}$ 
but starting from $Z^{\prime (p)}_{0}(0) =(b-1) \lfloor \ln^{4}n \rfloor + 1 - \Delta_{n}$. We recall that 
\begin{eqnarray*}
 W^{\prime}(t) := e^{-(b-1)t}Z^{\prime}(t) \hspace*{5mm} \text{and} \hspace*{5mm} W_{0}^{\prime (p)}(t) := e^{-(bp-1)t}Z^{\prime (p)}(t), \hspace*{5mm} t \geq 0
\end{eqnarray*}

\noindent are nonnegative square-integrable martingales which converge a.s. and in
$L^{2}(\mathbb{P})$. \\

\noindent {\bf Proof of  Lemma \ref{lema8}:} An application of Doob's inequality
(see, e.g., Equation (6) in \cite{Fort}) shows for all $\eta > 0$ that
\begin{eqnarray*} 
 \lim_{n \rightarrow \infty} \mathbb{P} \left( \left| e^{-(b-1) \tau^{\prime}( \ln^{4}n )}Z^{\prime}(\tau^{\prime}( \ln^{4}n ))- Z^{\prime}(0) \right| >  \eta \ln^{3} n  \right) = 0
\end{eqnarray*}

\noindent and using the fact that $Z_{0}^{\prime (p)}(0) \leq  (b-1)\lfloor \ln^{4} n \rfloor+1$, we also get
\begin{eqnarray*}
 \lim_{n \rightarrow \infty} \mathbb{P} \left( \left|  e^{-(bp-1) \tau^{\prime}( \ln^{4}n )}Z^{\prime (p)}_{0}(\tau^{\prime}( \ln^{4}n )) -Z_{0}^{\prime (p)}(0) \right| >  \eta \ln^{3} n  \right) = 0.
\end{eqnarray*}

Then, since $\Delta_{0,n} = Z^{\prime}(\tau^{\prime}( \ln^{4}n )) -(b-1)\lfloor \ln^{4} n \rfloor-1$, $\Delta_{n} = (b-1)\lfloor \ln^{4} n \rfloor+1 - Z^{\prime (p)}_{0}(0)$, and $Z^{\prime}(0) =(b-1) \lfloor \ln^{4} n \rfloor+1$, one readily gets 
\begin{eqnarray*}
 \Delta_{n}- \Delta_{0,n} =  Z^{\prime}(\tau^{\prime}( \ln^{4}n ))\left(e^{-(b-1) \tau^{\prime}( \ln^{4}n )} -1 \right) - Z^{\prime (p)}_{0}(\tau^{\prime}( \ln^{4}n )) \left(e^{-(bp-1) \tau^{\prime}( \ln^{4}n )} -1 \right) + o(\ln^{3}n)
\end{eqnarray*}

\noindent in probability. We next note from Lemma \ref{lema9} that 
\begin{eqnarray*}
 Z^{\prime (p)}_{0}(\tau^{\prime}( \ln^{4}n )) \left(e^{(1-p) \tau^{\prime}( \ln^{4}n )} -1 \right) = o(\ln^{3}n) \hspace*{5mm} \text{in probability},
\end{eqnarray*}

\noindent which yields
\begin{eqnarray*}
 \Delta_{n} - \Delta_{0,n}  =  \left( W^{\prime}(\tau^{\prime}( \ln^{4}n )) - W^{\prime (p)}_{0}(\tau^{\prime}( \ln^{4}n )) \right) \left(1-e^{(b-1) \tau^{\prime}( \ln^{4}n )} \right)  + o(\ln^{3}n) 
\end{eqnarray*}

\noindent in probability. Since by Lemma \ref{lema8} we have that 
\begin{eqnarray*}
\lim_{n \rightarrow \infty} \left(1-e^{(b-1) \tau^{\prime}( \ln^{4}n )} \right) = 0 \hspace*{6mm} \text{in probability},
\end{eqnarray*} 

\noindent we must verify that 
\begin{eqnarray*}
 W^{\prime}(\tau^{\prime}( \ln^{4}n )) - W^{\prime (p)}_{0}(\tau^{\prime}( \ln^{4}n )) = o(\ln^{3}n) \hspace*{6mm} \text{in probability},
\end{eqnarray*}

\noindent in order to get the result of Lemma \ref{lema8}. We observe from properties 
of square-integrable martingales that
\begin{eqnarray*}
\mathbb{E} \left[ \left( W^{\prime}(\tau^{\prime}( \ln^{4}n )) \right)^{2} \right] = \mathbb{E} \left[ \left[ W^{\prime} \right]_{\tau^{\prime}( \ln^{4}n )} \right] 
\end{eqnarray*}

\noindent where
\begin{eqnarray*}
\left[ W^{\prime} \right]_{t} = \sum_{0 \leq s \leq t} e^{-2(b-1)s}|Z^{\prime}(s)-Z^{\prime}(s-)|^{2} \hspace*{6mm} \text{for} \hspace*{3mm} t \geq 0. 
\end{eqnarray*}

A straightforward calculation shows that the compensator of jump process $\left[ W^{\prime} \right]$ is
\begin{eqnarray*}
\langle W^{\prime} \rangle_{t} = (b-1)^{2} \int_{0}^{t} e^{-2(b-1)s}Z^{\prime}(s){\rm d}s \hspace*{6mm} \text{for} \hspace*{3mm} t \geq 0, 
\end{eqnarray*}

\noindent that is $\left[ W^{\prime} \right]_{t} -\langle W^{\prime} \rangle_{t}$ is a local martingale. Thus,
\begin{eqnarray*}
\mathbb{E} \left[ \left( W^{\prime}(\tau^{\prime}( \ln^{4}n )) \right)^{2} \right] = \mathbb{E} \left[ \langle W^{\prime}  \rangle_{\tau^{\prime}( \ln^{4}n )} \right] \leq \mathbb{E} \left[ \langle W^{\prime}  \rangle_{\infty} \right] = (b-1)((b-1) \lfloor \ln^{4}n \rfloor +1) .
\end{eqnarray*}

\noindent Hence by the Markov inequality we have that
\begin{eqnarray*}
W^{\prime}(\tau^{\prime}( \ln^{4}n )) = o(\ln^{3} n) \hspace*{6mm} \text{in probability}.
\end{eqnarray*}

\noindent Similarly one gets
\begin{eqnarray*}
W^{\prime(p)}_{0}(\tau^{\prime}( \ln^{4}n )) = o(\ln^{3} n) \hspace*{6mm} \text{in probability},
\end{eqnarray*}

\noindent from where our claim follows. \hfill $\Box$ 

\subsection{The spread of fluctuations} \label{spread}

The purpose here is to resume the growth of the system of branching processes with rare mutation from the size 
$(b-1)\lfloor \ln^{4} n \rfloor+1$ to the size $(b-1) n + 1$ and observe that the germ of the fluctuations $\Delta_{n}$
spreads regularly. In this direction, we proceed similarly as the last part of the preceding section. 
We recall that $Z^{\prime}$ denotes the process of the total population
started from $Z^{\prime}(0) = (b-1)\lfloor \ln^{4} n \rfloor+1$. We consider
\begin{eqnarray*}
 \tau^{\prime}(n)= \inf \{t \geq 0: Z^{\prime}(t) = (b-1) n + 1 \},
\end{eqnarray*}

\noindent the time needed for the total population to reach size $(b-1) n + 1$. Hence, in the notation of Theorem \ref{teo2}
\begin{eqnarray*}
 G_{n} = Z_{0}^{\prime (p)}(\tau^{\prime}(n)),
\end{eqnarray*}

\noindent where as the previous section, we write $Z_{0}^{\prime (p)}$ for the process that counts
the number of individuals with the ancestral type starting from $Z_{0}^{\prime (p)}(0) 
= (b-1)\lfloor \ln^{4} n \rfloor+1-\Delta_{n}$. \\

We have now all the ingredients to establish Theorem \ref{teo2}. \\

\noindent {\bf Proof of Theorem \ref{teo2}:}  Again from the estimate of Equation (6) in \cite{Fort}, we get for all $\eta > 0$ that
\begin{eqnarray*}
 \lim_{n \rightarrow \infty} \mathbb{P} \left( \left| ((b-1) n + 1 )e^{-(b-1) \tau^{\prime}(n)}- (b-1)\ln^{4} n -1\right| >  \eta \ln^{3} n  \right) = 0,
\end{eqnarray*}

\noindent this yields
\begin{eqnarray}  \label{ec4}
e^{(b-1)\tau^{\prime}(n)} = \frac{n}{\ln^{4}n} + o \left( \frac{1}{\ln n} \right) \hspace*{6mm} \text{in probability.} 
\end{eqnarray}

On the other hand, using the fact that $Z_{0}^{\prime (p)}(0) \leq  (b-1)\lfloor \ln^{4} n \rfloor+1$, we also get
\begin{eqnarray*}
 \lim_{n \rightarrow \infty} \mathbb{P} \left( \left| e^{-(bp-1) \tau^{\prime}(n)}Z^{\prime (p)}_{0}(\tau^{\prime}(n))-Z_{0}^{\prime (p)}(0) \right| >  \eta \ln^{3} n  \right) = 0,
\end{eqnarray*}

\noindent and deduce that
\begin{eqnarray*}
 G_{n} = e^{(bp-1)\tau^{\prime}(n)}((b-1)\ln^{4}n - \Delta_{n})+  o \left(\frac{n}{\ln n} \right) \hspace*{5mm} \text{in probability.}
\end{eqnarray*}

Next, it is convenient to apply Skorokhod's representation theorem and assume that the weak convergence in Corollary \ref{cor2} holds almost surely. Hence
\begin{eqnarray*}
 G_{n} = e^{(bp-1)\tau^{\prime}(n)} \left((b-1)\ln^{4}n - \frac{\beta}{\beta -1} c \ln^{3}n \left(  3\ln \ln n  +  \left(\mathcal{Z}_{c,\beta} + 1 - \frac{1}{\beta}  \right)  \right) \right)+ o \left(\frac{n}{\ln n} \right)
\end{eqnarray*}

\noindent in probability. We next note from (\ref{ec4}) that
\begin{eqnarray*}
 e^{(bp-1)\tau^{\prime}(n)} =  e^{-\beta c}\frac{n}{\ln^{4}n}  +  4\beta c e^{-\beta c}n \frac{\ln \ln n}{\ln^{5}n} +o \left( \frac{n}{\ln^{5} n} \right)
\end{eqnarray*}

\noindent in probability. Recall that $\beta = b / (b-1)$, then
\begin{eqnarray*}
 G_{n}  =  \frac{1}{\beta-1}e^{-\beta c}n +   \frac{\beta}{\beta -1} ce^{-\beta c}n \frac{\ln \ln n}{\ln n} -  \frac{\beta}{\beta -1} c e^{-\beta c} \frac{n}{\ln n} \left( \mathcal{Z}_{c,\beta}  + 1 -  \frac{1}{\beta}  \right)  +  o \left(\frac{n}{\ln n} \right)
\end{eqnarray*}

\noindent in probability, which completes the proof. \hfill $\Box$ 

\section{Proof of Theorem \ref{teo1}} \label{sec3}

Our approach is based in the introduction of a continuous version of the construction of a $b$-ary recursive tree that enables us to superpose Bernoulli bond percolation dynamically in the tree structure. We begin at time $0$ from the tree with just one internal vertex which corresponds to the root having $b$ external vertices. Once the random tree  with size $n \geq 1$ has been constructed, we equip each of the $(b-1)n +1$  external vertices with independent exponential random variables $\zeta_{i}$ of parameter $1$. Then, after a waiting time equal to $\min_{i\in \{1, \dots, (b-1)n+1 \}} \zeta_{i}$, one of the external vertices is chosen uniformly at random and is replaced it by the internal vertex $n+1$ to which $b$ new leaves are attached. We observe that $\min_{i\in \{1, \dots, (b-1)n+1 \}} \zeta_{i}$ is exponentially distributed with  parameter $(b-1)n +1$. \\
 
We denote by $T(t)$ the tree which has been constructed at time $t \geq 0$, and by $|T(t)|$ its size, i.e. the number of internal vertices. The process of the size $(|T(t)|: t \geq 0)$ is clearly Markovian and if we define 
\begin{eqnarray*}
\gamma(n) = \inf \{ t \geq 0: |T(t)| = n\}, \hspace*{5mm} n \geq 1,
\end{eqnarray*}

\noindent then $T(\gamma(n))$ is a version of the $b$-ary recursive tree of size $n$, $T_{n}$. 
However for our purpose it will be more convenient work with the process $Y$ defined by
\begin{eqnarray} \label{ec20}
Y(t) = (b-1)|T(t)|+1, \hspace*{4mm} t \geq 0
\end{eqnarray}

\noindent with starting value $Y(0) = b$. It should be clear that $Y$ 
is a Yule process as described in Section \ref{sec2},
i.e. it has jumps of size $b-1$ and unit birth rate per unit population size.
We also point out that the process $Y$ gives us the number of external vertices on the 
tree. \\

We next superpose Bernoulli bond percolation with parameter $p=p_{n}$ defined in (\ref{ec1}) to the growth algorithm in continuous time  of the $b$-ary recursive tree. We follow the approach developed by Bertoin and Uribe Bravo \cite{Be1} but with a slight modification. We draw an independent Bernoulli random variable  $\epsilon_{p}$ with parameter $p$, each time an internal vertex is inserted. The edge which connects this new internal vertex is cut at its midpoint when $\epsilon_{p} = 0$ and remains intact otherwise. This disconnects the tree into connected clusters which motivates the following. We write $T^{(p)}(t)$ for the resulting combinatorial structure at time $t$. So, the percolation
clusters of $T(t)$ are the connected components by a path of intact edges of $T^{(p)}(t)$. \\

Let $T_{0}^{(p)}(t)$ be the subtree that contains the root.
We write $H_{0}^{(p)}(t)$ for the number of half-edges pertaining to the root cluster at time $t$. So that, its number of external vertices is given by
\begin{eqnarray*}
Y_{0}^{(p)}(t) = (b-1) |T_{0}^{(p)}(t)|+1-H_{0}^{(p)}(t).
\end{eqnarray*}


We are now be able to observe the connection with the system of branching processes 
with rare mutations described in the preceding section. It should be plain from the construction
that the size of the root-cluster at time $t$, i.e. $Y_{0}^{(p)}(t)$, of $T(t)$ after percolation with parameter
$p$, coincides with the number of individuals with the ancestral type $Z_{0}^{(p)}(t)$ in the system 
$\mathbf{Z}^{(p)}$ of branching processes with rare mutations of Section \ref{sec2}.
In fact, we already mentioned that the process $Y$ has the same random evolution as 
the process of the total size in the system $Z$. Recall that the algorithm for constructing a $b$-ary recursive tree is run until the time 
\begin{eqnarray*}
 \gamma(n) = \inf \{ t \geq 0: |T(t)| = n\} = \inf \{ t \geq 0: Y(t) = (b-1)n+1\} 
\end{eqnarray*}

\noindent when the structure has $n$ internal vertices. Then, the size $C_{0}^{(p)}$ of the percolation cluster 
containing the root when the tree has $n$ internal vertices satisfies
\begin{eqnarray*}
 C_{0}^{(p)} = |T_{0}^{(p)}(\gamma(n))|.
\end{eqnarray*}

In addition, it should be plain that 
\begin{eqnarray*}
 Y_{0}^{(p)}(\gamma(n)) = (b-1) C_{0}^{(p)} + 1 - H_{0}^{(p)}(\gamma(n)),
\end{eqnarray*}

\noindent coincides with the number of individuals with the ancestral type 
in the branching system $\mathbf{Z}^{(p)}$, at time when the total population
reaches the size $(b-1)n + 1$, i.e. $G_{n}$, according to the notation of Theorem \ref{teo2}.
Hence in order to establish Theorem \ref{teo1}, it is sufficient to get an estimate of the number 
of half-edges pertaining to the root-subtree at time $\gamma(n)$. 

\begin{lemma} \label{cor1}
 We have
 \begin{eqnarray*} 
 \lim_{n \rightarrow \infty} \frac{\ln n}{n} H^{(p)}_{0}(\gamma(n)) = ce^{-\beta c} \hspace*{5mm} \text{in probability.}
\end{eqnarray*}
\end{lemma}

\noindent {\bf Proof:}  We observe that the processes
\begin{eqnarray*}
 H_{0}^{(p)}(t) - (1-p) \int_{0}^{t} Y_{0}^{(p)}(s) {\rm d}s \hspace*{5mm} \text{and} \hspace*{5mm}  Y_{0}^{(p)}(t) - (bp-1) \int_{0}^{t} Y_{0}^{(p)}(s) {\rm d}s, \hspace*{5mm} t \geq 0
\end{eqnarray*}

\noindent are martingales. Thus,
\begin{eqnarray*}
 L^{(p)}(t) := H_{0}^{(p)}(t) - \frac{1-p}{bp-1} Y_{0}^{(p)}(t), \hspace*{5mm} t \geq 0
\end{eqnarray*}

\noindent is also a martingale. Observe that since $p = p_{n}$ satisfies (\ref{ec1}), for $n$ large enough such that $2/(b+1) \leq p \leq 1$, its jumps $|L^{(p)}(t) -L^{(p)}(t-)|$
have size at most $b$. Since there are at most $n$ jumps up to time $\gamma(n)$, the bracket of $L^{(p)}$ can be bounded by $[L^{(p)}]_{\gamma(n)} \leq b^{2}n$. Hence we have
\begin{eqnarray} \label{ec6}
 \lim_{n \rightarrow \infty} \mathbb{E} \left( \left| \frac{\ln  n}{n} L^{(p)}(\gamma(n))  \right|^{2} \right) = 0.
\end{eqnarray}

On the other hand, we know from Lemma \ref{lema4} that
\begin{eqnarray*}
 \lim_{n \rightarrow \infty} e^{-(b-1)\gamma(n)}Y(\gamma(n)) = \lim_{n \rightarrow \infty} e^{-(bp-1)\gamma(n)}Y_{0}^{(p)}(\gamma(n)) = W(\infty) \hspace*{5mm} \text{in probability}
\end{eqnarray*}

\noindent which implies that
\begin{eqnarray*} 
 \lim_{n \rightarrow \infty} \frac{Y_{0}^{(p)}(\gamma(n))}{n} = (b-1) e^{- \beta c} \hspace*{5mm} \text{in probability},
\end{eqnarray*}

\noindent and the result follows readily from (\ref{ec6}), the above limit and the fact that $1-p = o(1)$.  \hfill $\Box$ \\

Therefore, from the identity
\begin{eqnarray*}
 C_{0}^{(p)} = \frac{Y_{0}^{(p)}(\gamma(n))- 1 + H_{0}^{(p)}(\gamma(n))}{(b-1)},
\end{eqnarray*}

\noindent Theorem \ref{teo2} applies to $Y_{0}^{(p)}(\gamma(n))$ and Lemma \ref{cor1} yields the result of Theorem \ref{teo1}.

\section{Percolation on scale-free trees} \label{sec4}

We conclude this work by showing that the approach used in the proof 
of Theorem \ref{teo1} can be also applied to study percolation on scale-free trees, which form a family of random trees on a set of ordered vertices, say $\{0, 1, \dots, n \}$, that grow following a preferential attachment algorithm. Specifically, fix a parameter $a \in (-1, \infty)$, and start for $n=1$ from the tree $T_{1}^{(a)}$ on $\{0,1\}$ which has a single edge connecting $0$ and $1$. Suppose that  $T_{n}^{(a)}$ has been constructed for some $n \geq 2$, and for every $i \in \{ 0, 1, \dots, n\}$, denote by $d_{n}(i)$ the degree of the vertex $i$ in $T_{n}^{(a)}$. Then conditionally given $T_{n}^{(a)}$, the tree $T_{n+1}^{(a)}$ is built by adding an edge between the new vertex $n+1$ and a vertex $v_{n}$ in $T_{n}^{(a)}$ chosen at random according to the law
\begin{eqnarray*}
 \mathbb{P} \left( v_{n} = i| T_{n}^{(a)} \right) = \frac{d_{n}(i)+ a}{2n+a(n+1)}, \hspace*{5mm} i \in \{ 0, 1, \dots, n\}.
\end{eqnarray*}

Clearly, the preceding expression defines a probability measure since the sum of the degrees of a tree with $n+1$
vertices equals $2n$. Note also that when one lets $a\rightarrow \infty$ the algorithm yields
an uniform recursive tree since $v_{n}$ becomes uniformly distributed on $\{ 0, 1, \dots, n\}$. We then perform Bernoulli bond percolation with parameter given by (\ref{ec1}), i.e. $p_{n} = 1 - c /\ln n$, where $c>0$ is fixed. It has been observed by Bertoin and Uribe Bravo \cite{Be1} that this choice of the percolation parameter corresponds to the supercritical regime. More precisely, the size of the cluster $\Gamma^{(\alpha)}_{n}$ containing the root satisfies 
\begin{eqnarray*}
\lim_{n \rightarrow \infty} n^{-1} \Gamma^{(\alpha)}_{n} = e^{-\alpha c} \hspace*{6mm} \text{in probability},
\end{eqnarray*} 

\noindent where $\alpha = (1+a)/(2+a)$. We are interested in the fluctuations of $\Gamma^{(\alpha)}_{n}$, and show that an analogous result to Theorem \ref{teo1} holds for large scale-free random trees. 

\begin{theorem} \label{teo3}
  Set $\alpha = (1 + a)/(2 + a)$, and assume that the percolation parameter $p_{n}$ is given by (\ref{ec1}). Then as $n \rightarrow \infty$, there is the weak convergence
 \begin{eqnarray*}
  \left( n^{-1}\Gamma_{n}^{(\alpha)}- e^{- \alpha c} \right) \ln n - \alpha c e^{-\alpha c} \ln \ln n \Rightarrow  - \alpha c e^{-\alpha c}  \mathcal{Z}_{c,\alpha}^{\prime}
 \end{eqnarray*}
 \noindent where
\begin{eqnarray} \label{rec6}
\mathcal{Z}_{c,\alpha}^{\prime}=   \mathcal{Z} - \kappa_{\alpha}^{\prime} + \ln \left( \alpha c \right)  
\end{eqnarray}

\noindent with $\mathcal{Z}$ the continuous Luria-Delbr\"uck distribution and 
\begin{eqnarray*} 
\kappa_{\alpha}^{\prime} = 1- \frac{1}{\alpha} + \frac{1}{\alpha} \sum_{k=2}^{\infty} \frac{(\alpha)_{k}}{k!} \frac{(-1)^{k}}{k-1}.  
\end{eqnarray*}
\end{theorem}

We now focus on the proof of Theorem \ref{teo3}. We follow the route 
used in the proof of Theorem \ref{teo1}, and we analyze a system of branching processes with rare neutral mutations. 
We point out that in order to avoid repetitions, some technical details will be skipped. \\

We start by considering a pure birth branching process $Z^{(a)} = (Z^{(a)}(t):t \geq 0)$ in
continuous space, that has only jumps of size $2+a$, and with unit birth rate per unit population size. 
We shall be mainly interested in a class of population systems which arises by incorporating
mutations to the preceding branching process. More precisely, we describe the evolution of
such a system by a process $\mathbf{Z}^{(p, a)} = (\mathbf{Z}^{(p, a)} :t \geq 0)$,
where for each $t \geq 0$, $\mathbf{Z}^{(p, a)}(t) = (Z_{0}^{(p , a)}(t), Z_{1}^{(a)}(t), \dots)$
is a collection of nonnegative variables. At the initial time, all the sub-populations $Z_{i}^{(a)}(0)$ of type 
$i \geq 1$ are taken to be equal to zero, and $Z_{0}^{(p , a)}(0) = 2 + 2a$ which is the 
size of the ancestral (or clone) population. We consider that at rate $p$ per unit population
size, the clone population produces $2+a$ new clones, and that at rate $1-p$ per unit population size, 
they always create a single mutant population of a new type of size $1+a$. The new mutant populations behave 
as the process $Z^{(a)}$ but starting from $1+a$. Clearly, the sum over all sub-populations 
\begin{eqnarray*}
 Z^{(a)}(t) = Z_{0}^{(p, a)}(t) + \sum_{i \geq 1} Z_{i}^{(a)}(t), \hspace*{6mm} t \geq 0,
\end{eqnarray*}

\noindent evolves as the pure birth branching process described at the beginning of this paragraph. \\

We next observe the growth of the system of branching process $\mathbf{Z}^{(p, a)}$ until the time
\begin{eqnarray*}
 \tau^{(a)}(\ln^{4}n) = \inf\{ t \geq 0: Z^{(a)}(t) = (2+a)\lfloor \ln^{4}n \rfloor + a \},
\end{eqnarray*}

\noindent which is when the total size of the population reaches $(2+a)\lfloor \ln^{4}n \rfloor + a$. Our first purpose is to estimate precisely the number $\Delta_{n}^{(\alpha)}$ of mutants at this time. This stage corresponds to the analysis of the germ, and approach line will be similar to that in Section \ref{germ}. In this direction, it will be useful to study the number of mutants $\Delta_{0,n}^{(\alpha)}$ at time
\begin{eqnarray*}
 \tau^{(a)}_{0}(\ln^{4}n) = \inf\{ t \geq 0: Z^{(p,a)}_{0}(t) = (2+a)\lfloor \ln^{4}n \rfloor + a \}
\end{eqnarray*}

\noindent whose distribution is easier to estimate than that of $\Delta_{n}^{(\alpha)}$. 
We shall establish the following limit theorem in law that is equivalent to the Proposition
\ref{pro1}.

\begin{proposition} \label{rpro1}
 As $n \rightarrow \infty$, there is the weak convergence
 \begin{eqnarray*}
  \frac{\Delta_{0,n}^{(\alpha)}}{\ln^{3}n} - 3 \frac{\alpha}{1-\alpha} c \ln \ln n \Rightarrow \frac{\alpha}{1-\alpha} c \left( \mathcal{Z}_{c, \alpha}^{\prime} + 1-\frac{1}{\alpha}\right)
 \end{eqnarray*}
\noindent where $\mathcal{Z}_{c, \alpha}^{\prime}$ is the random variable defined in (\ref{rec6}).
\end{proposition}

Before proving the Proposition \ref{rpro1}, it is convenient to introduce the following
representation of the total mutant population as we have done in Section \ref{germ}. For $i \geq 1$,
we write 
\begin{eqnarray*}
 b_{i}^{(p)} = \inf \{t \geq 0: Z_{i}^{(a)}(t) > 0 \},
\end{eqnarray*}

\noindent for the birth time of the sub-population with type $i$. Then the processes 
$(Z_{i}^{(a)}(t - b_{i}^{(p)}): t \geq b_{i}^{(p)})$ form a sequence of i.i.d. 
branching processes with the same law as $Z^{(a)}$ but starting value $1+a$, which 
is independent of the birth-times $(b_{i}^{(p)})_{i \geq 1}$ and the process $Z_{0}^{(p, a)}$.
Moreover, this sequence is also independent of the process that counts the number of mutation
events which is defined by $M^{(a)}(t) = \max \{ i \geq 1: Z_{i}^{(a)}(t) > 0 \}$. Thus, since
the jump times of $M^{(a)}$ are in fact $b_{1}^{(p)} < b_{2}^{(p)} < \dots$, we can express the total
mutant population at time $t \geq 0$ as,
\begin{eqnarray*}
 Z^{(a)}_{\text{m}}(t) = \sum_{i=1}^{M^{(a)}(t)} Z_{i}^{(a)}(t - b_{i}^{(p)}).
\end{eqnarray*}

We observe that for $i \geq 1$, the process $((2+a)^{-1}Z_{i}^{(a)}(t - b_{i}^{(p)}): t \geq b_{i}^{(p)})$
is a Yule branching process in continuous space with birth rate $2+a$ per unit population size.
Then similarly as we obtained the result in Lemma \ref{lema7}, we get for $t\geq 0$ and $\theta \in \mathbb{R}$,
\begin{eqnarray} \label{rcar}
 \mathbb{E}[e^{i \theta Z^{(a)}_{\text{m}}(t)}] = \mathbb{E} \left[ \exp \left((1-p) \int_{0}^{t} Z_{0}^{(p, a)}(t-s) (\varphi_{s}^{(a)}(\theta)-1) {\rm d}s\right)\right]
\end{eqnarray}

\noindent where 
\begin{eqnarray*} 
 \varphi_{t}^{(a)}(\theta) = \mathbb{E}\left[ e^{i \theta Z^{(a)}(t)} \big| Z^{(a)}(0) = 1+a \right] = \left( \frac{e^{i \theta(2+a)} e^{-(2+a)t}}{  1 - e^{i \theta(2+a)} + e^{i \theta(2+a)} e^{-(2+a)t}} \right)^{\alpha}
\end{eqnarray*} 

\noindent with $\alpha = (1+a)/(2+a)$. \\

At this point, the difference between the constants in Theorem \ref{teo1} and 
Theorem \ref{teo3} must be evident, mostly due to the different behavior of the 
branching processes associated to the $b$-ary recursive trees and scale-free 
random trees. Essentially, the constant $\kappa_\beta$ of Theorem \ref{teo1}
depends of the characteristic function (\ref{ec11}) through the computations made in 
Lemma \ref{lema5}, which is clearly distinct from (\ref{rcar}). \\

We are now able to establish Proposition \ref{rpro1}. \\

\noindent {\bf Proof of Proposition \ref{rpro1}:} 
We fix $\theta \in \mathbb{R}$ and define $m_{n} = \alpha c \ln^{3} n$.
Since we have the identity $\Delta_{0,n}^{(\alpha)} = Z^{(a)}_{\text{m}}(\tau_{0}^{(a)} (\ln^{4} n))$,
it follows from (\ref{rcar}) that the characteristic function of $m^{-1}_{n} \Delta_{0,n}^{(\alpha)} - 
(1-\alpha)^{-1} \ln m_{n}$ is given by
\begin{eqnarray} \label{dec2}
 \mathbb{E}\left[e^{i \theta ( m^{-1}_{n} \Delta_{0,n}^{(\alpha)} - 
(1-\alpha)^{-1} \ln m_{n}) } \right] = \mathbb{E}\left[ \exp \left( I^{(p, a)}(\tau_{0}^{(a)} (\ln^{4}n) ) - (1-\alpha)^{-1} \ln m_{n} \right)\right],
\end{eqnarray}

\noindent where 
\begin{eqnarray*}
 I^{(p, a)}(t) = (1-p) \int_{0}^{t} Z_{0}^{(p, a)}(t-s) (\varphi_{s}^{(a)}(u) -1) {\rm d}s \hspace*{5mm} \text{for} \hspace*{2mm} t \geq 0,
\end{eqnarray*}

\noindent and $u = \theta / (\alpha c \ln^{3}n)$. Next, a similar computation as the proof of Lemma 
\ref{lema1} shows that
\begin{eqnarray} \label{dec1}
 \lim_{n \rightarrow \infty} \left(I^{(p, a)}(\tau_{0}^{(a)} (\ln^{4} n) ) - I^{(p, a)}_{{\rm m}}(\tau_{0}^{(a)} (\ln^{4} n) ) \right) = 0 \hspace*{6mm} \text{in probability}, 
\end{eqnarray}

\noindent where 
\begin{eqnarray*}
 I^{(p, a)}_{{\rm m}}(t) = (1-p) W_{0}^{(p, a)}(\infty) e^{(2+a)t}\int_{0}^{t} e^{-(2+a)s} (\varphi_{s}^{(a)}(u) -1) {\rm d}s \hspace*{5mm} \text{for} \hspace*{2mm} t \geq 0.
\end{eqnarray*}

\noindent and  $W_{0}^{(p, a)}(\infty)$ is defined as the terminal value of the martingale $W_{0}^{(p, a)}(t) = e^{-(1+p(1+a))t}Z_{0}^{(p, a)}(t)$. Moreover, the integral of the previous expression can be computed explicitly,
\begin{eqnarray*}
 \int_{0}^{t} e^{-(2+a)s} (\varphi_{s}^{(a)}(u) -1) {\rm d}s = \frac{1-e^{iu(2+a)}}{(2+a)e^{iu(2+a)}} \left( \alpha \ln(1 - e^{iu(2+a)} + e^{iu(2+a)} e^{-(2+a)t}) + \kappa_{\alpha,u}^{\prime}(t)\right),
\end{eqnarray*}

\noindent with
\begin{eqnarray*}
 \kappa_{\alpha,u}^{\prime}(t) = \sum_{k=2}^{\infty} \frac{\left(\alpha\right)_{k}}{k!}   \frac{(e^{iu(2+a)}-1)^{k-1}}{k-1} \left( 1 - \frac{1}{(1-e^{iu(2+a)}+e^{iu(2+a)}e^{-(2+a)t})^{k-1}} \right).
\end{eqnarray*}

\noindent Hence from Lemma 3 in \cite{Be1} (which is the analog of Lemma \ref{lema4}), 
we conclude after some computations that
\begin{align*}
  & \lim_{n \rightarrow \infty} \left(I^{(p, a)}_{\text{m}} (\tau_{0}^{(a)}(\ln^{4} n )) - i\theta (1-\alpha)^{-1} \ln m_{n} \right) \\
  &~~~~~~ = -  i\theta (1-\alpha)^{-1} \left( \kappa_{\alpha}^{\prime} -1 + \frac{1}{\alpha} \right)  - i \theta(1-\alpha)^{-1} \ln |(1-\alpha)^{-1}\theta| - \frac{1}{2} \pi |(1-\alpha)^{-1}\theta| 
\end{align*} 

\noindent in probability, and our claim follows from (\ref{dec1}), by letting $n \rightarrow \infty$ in (\ref{dec2}). \hfill $\Box$ \\

It follows now readily from the same arguments that we have developed to show Lemma \ref{lema8} that $\Delta_{n}^{(\alpha)}$ and $\Delta_{0,n}^{(\alpha)}$ have the same asymptotic
behavior. Specifically, we have:
\begin{corollary} \label{rcor1}
 As $n \rightarrow \infty$, there is the weak convergence
 \begin{eqnarray*}
  \frac{\Delta_{n}^{(\alpha)}}{\ln^{3}n} - 3 \frac{\alpha}{1-\alpha} c \ln \ln n \Rightarrow \frac{\alpha}{1-\alpha} c \left( \mathcal{Z}_{c, \alpha}^{\prime} + 1-\frac{1}{\alpha}\right)
 \end{eqnarray*}
\noindent where $\mathcal{Z}_{c, \alpha}^{\prime}$ is the random variable defined in (\ref{rec6}).
\end{corollary}

Similarly as in Section \ref{spread}, we now resume the growth of the system of branching process with rare mutation from the size $(2+a)\lfloor \ln^{4}n \rfloor + a$ to the size $(2+a)n + a$, and show that the fluctuations of $\Delta_{n}^{(\alpha)}$ spread regularly. In this direction, we write
$Z^{\prime (a)} = (Z^{ \prime (a)}(t): t \geq 0)$ for the process
of the total size of the population started from $Z^{ \prime (a)}(0) = (2+a)\lfloor \ln^{4}n \rfloor + a$, which has the same
law as the branching process $Z^{(a)}$. We observe that the population at the time when we restart our observation consists of $\Delta_{n}^{(\alpha)}$ mutants and 
$(2+a)\lfloor \ln^{4}n \rfloor + a - \Delta_{n}^{(\alpha)}$
individuals with the ancestral type. Then, we write $Z^{\prime (p, a)}_{0} = (Z^{\prime (p, a)}_{0}(t): t \geq 0)$ for the process
that counts the number of clone individuals, which has the same law as $Z_{0}^{(p, a)}$ 
but starting from $Z^{\prime (p, a)}_{0}(0) =(2+a)\lfloor \ln^{4}n \rfloor + a - \Delta_{n}^{(\alpha)}$. We consider the time
\begin{eqnarray*}
 \tau^{\prime(a)}(n ) = \inf \{t \geq 0: Z^{\prime (a)}(t) = (2+a)n + a\}.
\end{eqnarray*}

Then the number of individuals with the ancestral type at time when the total population
generated by the branching process reaches $(2+a)n + a$ is given by
\begin{eqnarray*}
 G_{n}^{(\alpha)} = Z^{\prime (p,a)}_{0}(\tau^{\prime(a)}(n )).
\end{eqnarray*}

We are now able to state the following analog of Theorem \ref{teo2}.

\begin{theorem} \label{rteo2}
 Set $\alpha = (1+a)/(2+a)$. As $n \rightarrow \infty$, there is the weak convergence
 \begin{eqnarray*}
  \left(n^{-1}G_{n}^{(\alpha)}-\frac{1}{1-\alpha} e^{-\alpha c} \right) \ln n - \frac{\alpha}{1-\alpha} c e^{-\alpha c} \ln \ln n \Rightarrow  -  \frac{\alpha}{1-\alpha} c e^{-\alpha c} \left( \mathcal{Z}_{c,\alpha}^{\prime} +1- \frac{1}{\alpha}  \right),
 \end{eqnarray*}
 
 \noindent where $\mathcal{Z}_{c,\alpha}^{\prime}$ is the random variable defined in (\ref{rec6}).
\end{theorem}

\noindent {\bf Proof:} We recall that 
\begin{eqnarray*}
 W^{\prime(a)}(t) := e^{-(2+a)t}Z^{\prime(a)}(t) \hspace*{5mm} \text{and} \hspace*{5mm} W_{0}^{\prime (p, a)}(t) := e^{-(1+p(1+a))t}Z^{\prime (p, a)}(t), \hspace*{5mm} t \geq 0
\end{eqnarray*}

\noindent are nonnegative square-integrable martingales which converge a.s. and in
$L^{2}(\mathbb{P})$. Hence from the estimate of Equation (6) in \cite{Fort}, we get for all $\eta > 0$ that
\begin{eqnarray*}
 \lim_{n \rightarrow \infty} \mathbb{P} \left( \left| ((2+a)n+ a)e^{-(2+a) \tau^{\prime(a)}(n)}- ((2+a) \lfloor \ln^{4}n \rfloor +a)\right| >  \eta \ln^{3} n  \right) = 0,
\end{eqnarray*}

\noindent this yields
\begin{eqnarray*}  
e^{(2+a)\tau^{\prime (a)}(n)} = \frac{n}{\ln^{4}n} + o \left( \frac{1}{\ln n} \right) \hspace*{6mm} \text{in probability.} 
\end{eqnarray*}

On the other hand, using the fact that $Z_{0}^{\prime (p, a)}(0) \leq  (2+a)\lfloor \ln^{4}n \rfloor + a$, we also get
\begin{eqnarray*}
 \lim_{n \rightarrow \infty} \mathbb{P} \left( \left| e^{-(1+p(1+a)) \tau^{\prime(a)}(n)}Z^{\prime (p, a)}_{0}(\tau^{\prime (a)}(n))-Z_{0}^{\prime (p, a)}(0) \right| >  \eta \ln^{3} n  \right) = 0,
\end{eqnarray*}

\noindent and deduce that
\begin{eqnarray*}
 G_{n}^{(\alpha)} = e^{(1+p(1+a))\tau^{\prime(a)}(n)}((2+a) \ln^{4}n - \Delta_{n}^{(\alpha)})+  o \left(\frac{n}{\ln n} \right) \hspace*{5mm} \text{in probability.}
\end{eqnarray*}

\noindent Skorokhod's representation theorem enables us to assume that the weak convergence in Corollary \ref{rcor1} holds almost surely. Hence
\begin{eqnarray*}
 G_{n}^{(\alpha)} = e^{(1+p(1+a))\tau^{\prime(a)}(n)} \left((2+a)\ln^{4}n - \frac{\alpha}{1 - \alpha} c \ln^{3}n \left(  3\ln \ln n  +  \left(\mathcal{Z}_{c,\alpha}^{\prime}+1-\frac{1}{\alpha}  \right)  \right) \right)+ o \left(\frac{n}{\ln n} \right)
\end{eqnarray*}

\noindent in probability. It follows that
\begin{eqnarray*}
 G_{n}^{(\alpha)}  =  \frac{1}{1-\alpha}e^{-\alpha c}n +   \frac{\alpha}{1 - \alpha} ce^{-\alpha c}n \frac{\ln \ln n}{\ln n} -  \frac{\alpha}{1 - \alpha} c e^{-\alpha c} \frac{n}{\ln n} \left( \mathcal{Z}_{c,\alpha}^{\prime}  +1-\frac{1}{\alpha} \right)  +  o \left(\frac{n}{\ln n} \right)
\end{eqnarray*}

\noindent in probability, which completes the proof. \hfill $\Box$ \\

We have now all the ingredients to establish Theorem \ref{teo3}.  \\

\noindent {\bf Proof of Theorem \ref{teo3}:} We follow Bertoin and Uribe Bravo \cite{Be1},
and we consider a continuous time version of the growth 
algorithm with preferential attachment as we have done for the $b$-ary recursive trees. 
We start at $0$ from the tree $\{0,1\}$, and once the random tree with size
$n \geq 2$ has been constructed, we equip each vertex $i \in \{0, 1, \dots , n\}$ with 
an exponential random variable $\zeta_{i}$ of parameter $d_{n}(i) + a$, independently
of the other vertices. Then the next vertex $n + 1$ is attached after time $\min_{i \in \{0, 1, \dots, n\}} \zeta_i$ at
the vertex $v_{n} = \text{argmin}_{i \in \{0,1, \dots,n\}} \zeta_i$. Let us denote by 
$T^{(a)}(t)$ the tree which has been constructed at time $t$, and by $|T^{(a)}(t)|$ its size,
i.e. its number of vertices. It should be plain that if we define 
\begin{eqnarray*}
\gamma^{(a)}(n) = \inf \{t \geq 0: |T^{(a)}(t)| = n +1 \},
\end{eqnarray*}

\noindent then $T^{(a)}\left(\gamma^{(a)}(n)\right)$ is a version of a scale-free tree of size $n+1$, $T_{n}^{(a)}$.
Furthermore, the process $Y^{(a)}$ defined by
\begin{eqnarray*}
 Y^{(a)}(t) = (2 + a) |T^{(a)}(t)| -2, \hspace*{5mm} t \geq 0,
\end{eqnarray*}

\noindent is a pure branching process with initial value $Y^{(a)}(0) = 2a + 2$ that has only jumps of size $2 + a$, and
with unit birth rate per unit population size. Then we incorporate Bernoulli bond 
percolation to the algorithm similarly to how we did in Section \ref{sec2} for
the $b$-ary recursive trees. We draw an independent Bernoulli random variable $\epsilon_{p}$
with parameter $p$, each time an edge is inserted. If $\epsilon_{p} = 1$, the edge is left
intact, otherwise we cut this edge at its midpoint. We write $T^{(p,a)}(t)$ for the resulting combinatorial structure at time $t$.
Hence the percolation clusters of $T^{(a)}(t)$ are the connected components by a path of intact edges of $T^{(p,a)}(t)$. Let $T_{0}^{(p, a)}(t)$ 
be the subtree that contains the root. We write $H_{0}^{(p, a)}(t)$ for the number of half-edges pertaining to the root cluster at time $t$
and set 
\begin{eqnarray*}
Y_{0}^{(p, a)}(t) = (2+a) |T_{0}^{(p, a)}(t)| + H_{0}^{(p, a)}(t) - 2.
\end{eqnarray*}

We now observe the connection with the system of branching processes 
with rare mutations $\mathbf{Z}^{(p, a)}$. It should be plain from the construction
that $Y_{0}^{(p, a)}$ has the same random evolution as the process of the number of individuals with the ancestral type 
$Z_{0}^{(p,a)}(t)$. In fact, the process $Y^{(a)}$ coincides with the process of the total size $Z^{(a)}$. 
Then, the size $\Gamma_{n}^{(\alpha)}$ of the percolation cluster 
containing the root when the structure has size $n+1$ satisfies
$\Gamma_{n}^{(\alpha)} = |T_{0}^{(p, a)}(\gamma^{(a)}(n))|$. In addition, it should be plain that 
\begin{eqnarray} \label{rid}
 Y_{0}^{(p, a)}(\gamma^{(a)}(n)) = (2+a)\Gamma_{n}^{(\alpha)}  + H_{0}^{(p, a)}(\gamma^{(a)}(n)) - 2,
\end{eqnarray}

\noindent coincides with the number of individuals with the ancestral type 
in the branching system $\mathbf{Z}^{(p, a)}$, at time when the total population
reaches the size $(2+a)n + a$, i.e. $G_{n}^{(\alpha)}$. Then, in order to establish 
Theorem \ref{teo3}, it is sufficient to get an estimate of the number 
of half-edges pertaining to the root-subtree at time $\gamma^{(a)}(n)$. 
We follow the route of Lemma \ref{cor1} and observe that the process
\begin{eqnarray*}
 L^{(p,a)}(t) := H_{0}^{(p,a)}(t) - \frac{1-p}{1+p+pa} Y_{0}^{(p,a)}(t), \hspace*{5mm} t \geq 0
\end{eqnarray*}

\noindent is a martingale whose jumps have size at most $2+a$. Since there are at most 
$n$ jumps up to time $\gamma^{(a)}(n)$, the bracket of $L^{(p,a)}$ can be bounded by $[L^{(p,a)}]_{\gamma^{(a)}(n)} \leq (2+a)^{2}n$. Hence we have
\begin{eqnarray*} 
 \lim_{n \rightarrow \infty} \mathbb{E} \left( \left| \frac{\ln  n}{n} L^{(p,a)}(\gamma^{(a)}(n))  \right|^{2} \right) = 0.
\end{eqnarray*}

On the other hand, from Lemma 3 in \cite{Be1} we get that
\begin{eqnarray*}
 \lim_{n \rightarrow \infty} e^{-(2+a)\gamma^{(a)}(n)}Y^{(a)}(\gamma^{(a)}(n)) = \lim_{n \rightarrow \infty} e^{-(1+p(1+a))\gamma^{(a)}(n)}Y_{0}^{(p,a)}(\gamma^{(a)}(n)) = W^{(a)}(\infty) 
\end{eqnarray*}

\noindent in probability, where $W^{(a)}(\infty)$ is the terminal value of the martingale $W^{(a)}(t) = e^{-(2+a)t}Y^{(a)}(t)$. Thus, we have that
\begin{eqnarray*} 
 \lim_{n \rightarrow \infty} \frac{\ln n}{n} H_{0}^{(p,a)}(\gamma^{(a)}(n)) = c e^{- \alpha c} \hspace*{5mm} \text{in probability},
\end{eqnarray*}

\noindent and the result in Theorem \ref{teo3} follows from Theorem \ref{rteo2}
and the identity (\ref{rid}). \hfill $\Box$

\paragraph{Acknowledgements.} I would like to thank Jean Bertoin, who introduced me to 
the problem and gave guidance throughout the work. I would also like to 
thank the referee for 
his/her very careful review of the first draft of this work and useful 
comments which helped improve this manuscript.
 
 \bibliography{Gab}

\begin{thebibliography}{10}

\bibitem{Athe}
Athreya, K.~B. and Ney, P.~E.
\newblock {\em Branching processes}.
\newblock Dover Publications, Inc., Mineola, NY, 2004.
\newblock Reprint of the 1972 original [Springer, New York; MR0373040].

\bibitem{Be3}
Bertoin, J.
\newblock Almost giant clusters for percolation on large trees with logarithmic
  heights.
\newblock {\em J. Appl. Probab.}, 50(3):603--611, 2013.

\bibitem{Be4}
Bertoin, J.
\newblock On the non-gaussian fluctuations of the giant cluster for percolation
  on random recursive trees.
\newblock {\em Electron. J. Probab.}, 19:no. 24, 1--15, 2014.

\bibitem{Be2}
Bertoin, J.
\newblock Sizes of the largest clusters for supercritical percolation on random
  recursive trees.
\newblock {\em Random Structures Algorithms}, 44(1):29--44, 2014.

\bibitem{Be1}
Bertoin, J. and Uribe~Bravo, G.
\newblock Supercritical percolation on large scale-free random trees.
\newblock {\em Ann. Appl. Probab.}, 25(1):81--103, 2015.

\bibitem{Chauvin}
Chauvin, B., Klein, T., Marckert, J.-F., and Rouault, A.
\newblock Martingales and profile of binary search trees.
\newblock {\em Electron. J. Probab.}, 10:no. 12, 420--435, 2005.

\bibitem{Fort}
de~La~Fortelle, A.
\newblock Yule process sample path asymptotics.
\newblock {\em Electron. Comm. Probab.}, 11:193--199, 2006.

\bibitem{Drmota}
Drmota, M.
\newblock {\em Random trees}.
\newblock Springer. New York, Vienna, 2009.

\bibitem{Go}
Goldschmidt, C. and Martin, J.~B.
\newblock Random recursive trees and the {B}olthausen-{S}znitman coalescent.
\newblock {\em Electron. J. Probab.}, 10:no. 21, 718--745, 2005.

\bibitem{Hol}
Holmgren, C.
\newblock Random records and cuttings in binary search trees.
\newblock {\em Combin. Probab. Comput.}, 19(3):391--424, 2010.

\bibitem{Iks}
Iksanov, A. and M{\"o}hle, M.
\newblock A probabilistic proof of a weak limit law for the number of cuts
  needed to isolate the root of a random recursive tree.
\newblock {\em Electron. Comm. Probab.}, 12:28--35, 2007.

\bibitem{java}
Javanian, M. and Vahidi-Asl, M.~Q.
\newblock Depth of nodes in random recursive {$k$}-ary trees.
\newblock {\em Inform. Process. Lett.}, 98(3):115--118, 2006.

\bibitem{Ma}
Mahmoud, H.~M.
\newblock {\em Evolution of random search trees}.
\newblock John Wiley, New York, 1992.

\bibitem{parzen}
Parzen, E.
\newblock {\em Stochastic processes}.
\newblock Society for Industrial and Applied Mathematics (SIAM), Philadelphia,
  PA, 1999.
\newblock Reprint of the 1962 original.

\bibitem{Sch}
Schweinsberg, J.
\newblock Dynamics of the evolving {B}olthausen-{S}znitman coalescent.
\newblock {\em Electron. J. Probab.}, 17:no. 91, 50, 2012.

\end{thebibliography}
 \bibliographystyle{Myplain}

\end{document}